\UseRawInputEncoding
\documentclass[12pt,leqno]{article}
\usepackage{latexsym}
\usepackage{enumerate,amsmath,amsthm,amssymb,amsfonts,amscd}

\newtheorem{defn}{Definition}[section]
\newtheorem{thm}[defn]{Theorem}
\newtheorem{cor}[defn]{Corollary}
\newtheorem{prop}[defn]{Proposition}
\newtheorem{lem}[defn]{Lemma}

\numberwithin{equation}{section}

\def\F{\mbox{${\cal F}$}}

\newcommand{\cc}{\mathbb C}
\newcommand{\zz}{\mathbb Z}

\renewcommand{\a}{\alpha}

\newcommand{\ra}{\rightarrow}

\input{arrow.tex}

\begin{document}
\title{\bf Equivalency of the Corona problem and Gleason's problem in the theory of SCV}
\author{S. R. PATEL}
\maketitle
\begin{tabbing}
\> Secondary: 32A17, 32A38\\
\> 46E25, 46J10, 46J20
\end{tabbing}
\newpage
\noindent {\bf Abstract.} We establish an equivalency of the
Corona problem (1962) and Gleason's problem (1964) in the theory
of several complex variables. As an application, we give an
affirmative solution of the Corona problem for certain bounded
pseudoconvex domains or polydomains in $\cc^n$ including balls and
polydiscs. Indeed, we extend our recent work on Gleason's problem
based on the functional analytic approach, as well as extend
recent work of Clos.

We also use this equivalency or else other (functional analytic)
methods to affirmatively solve both problems for various Banach
spaces of bounded holomorphic functions (including certain
holomorphic mixed-norm spaces) on various types of domains in
$\cc^n$ such as holomorphic H\"{o}lder and Lipschitz spaces (left
open by Forn\ae ss and \O vrelid in 1983), holomorphic mean
Besov-Lipschitz spaces, Besov-Lipschitz spaces, Hardy-Sobolev
spaces, and a weighted Bergman space. The discussion goes via
first studying {\it Lipschitz algebras of holomorphic functions of
order $\a$}, where $\a\,\in\,(0,\,1]$; in particular, the Gel'fand
theory and the maximal ideal spaces of these algebras are
discussed.
\newpage
\noindent {\bf Structure of the article, with comments on the
strategies.} The organization of the paper is as follows. In \S 1,
we give an introduction to the corona problem and Gleason's
problem from the theory of several complex variables (briefly:
SCV), in order to formulate an equivalency by giving necessary and
sufficient background for this motive. In \S 2, we discuss an
equivalency between these two problems in Thms. 2.3, 2.4 and 2.6.
This has become possible as some results from [P4] are the base of
the idea of an equivalency, as well as certain results of Clos,
and Kerzmann and Nagel are extended. Then, as an application of
this equivalency or by applying ceratin functional analytic
methods such as the Banach algebra method as well as the classical
dead-end space method, in \S 3, we obtain the Gleason solvability
for various familiar Banach spaces of holomorphic functions,
including the affirmative solution of Forn\ae ss and \O vrelid's
question left open in 1983 for holomorphic H\"{o}lder and
Lipschitz spaces; Lem. 3.2, Thms. 3.13 and 3.14 are important
results in this view. Mostly, our results are complemented in the
sense that we do not recover the old results, obtained for other
known Banach spaces of holomorphic functions.

In \S 4, we solve the original corona problem by applying a
mixture of our Banach algebra method as well as the most elegant
method of Krantz and Li [KL1]. Also, we apply the same methods to
solve the corona problem for the same familiar Banach spaces of
bounded holomorphic functions on various types of domains in
$\cc^n$ (including the ball and the polydisc). In the process, we
also obtain extensions of certain results in some cases (e.g., see
[KL2, K5]). In \S 5, we give some concluding remarks on our
approaches, leading to future directions/open problems. Throughout
the paper, ``algebra" will mean a complex, commutative algebra
with identity.
\section{History of the two problems and formulation of an equivalency.}
\noindent {\bf Corona problem.} In 1938, Gel'fand discussed Banach
algebras in his thesis. In 1942, Kakutani asked the following
fundamental question (known as {\it Corona problem}) in the
subject: Let $H^\infty(D)$ be the Banach algebra of all bounded,
holomorphic functions on the open unit disc $D$ in the complex
plane $\cc$, equipped with the usual sup-norm, and let $f_j
\,\in\,H^\infty(D), \,j = 1, 2,\dots, n$, with the property that
$\sum_{j = 1}^n|f_j(z)|\,>\,\epsilon$ for some positive constant
$\epsilon$ and all $z\,\in\,D$ (call it a finite set of corona
data on the disc $D$), do there exist $g_j \,\in\,H^\infty(D), \,j
= 1, 2,\dots, n$, such that $\sum_{j = 1}^nf_jg_j\,=\,1$ (i.e.,
$\sum_{j = 1}^nf_j(z)g_j(z)\,=\,1$ for all $z\,\in\,D$)? This
latter equation is known as the B\'{e}zout equation [Ho, Chap. X]
and [G, Thm. V.1.8]. Recall that the B\'{e}zout's identity and
theorem have roots in arithmetic and abstract algebra. In fact,
this identity and the fundamental theorem of algebra imply the
result similar to an answer of the above fundamental question, but
for any number of polynomials in $n$ variables, known as Hilbert's
Nullstellensatz, establishing a fundamental relationship between
algebra and geometry -- the relationship being the basis of
algebraic geometry, and thus, they also have roots in algebraic
geometry.

The equivalent formulation of this question is to ask whether $D$
is dense in the maximal ideal space (aka spectrum)
$M(H^\infty(D))$ w.r.t. the Gel'fand (i.e., relative weak*-)
topology (equivalently, whether the corona is empty). This
equivalence is explained in some detail in [K4]. An authoritative
treatment of the corona problem on the disc appears in [G].
Carleson affirmatively solved the problem on the disc $D$ in 1962
[C]. By a domain we mean an open, connected and bounded set $G$ in
$\cc^n$. At this time there is no known domain in the plane $\cc$
on which the corona problem is failed.

Since 1962, the corona problem is ``automatically" open for the
balls and the polydiscs in $\cc^n,\,n\,>\,1$. The case for
arbitrary domains in $\cc^n,\,n \,>\,1$ is a lot more complicated,
and there have been investigations of this problem on various
types of domains; e.g., see some partial results on the corona
problem in SCV in [BCL, KL1-3, K5]. The extensive exposition of
history of the research of the problem is outlined in [DKSTW] --
in this volume, it is notable to recall Brudnyi's results related
to this problem for $H^\infty(R)$, where $R$ is a Caratheodory
hyperbolic Riemann surface [Br].

Recently, the corona problem has been solved for some basic
domains in $\cc^n$ such as $L^\infty$-pseudoconvex domains [T]
(recall that, for $n = 2$, this result was obtained by Krantz
[K5], and in [Pr], a similar conclusion is obtained under a much
more restrictive assumption that there exist $j\,\neq\,k$ such
that $\overline{|f_j|^{-1}(-\epsilon, \epsilon)} \bigcap
\overline{|f_k|^{-1}(-\epsilon, \epsilon)} \bigcap \partial G =
\emptyset$, where $\partial G$ denotes the boundary).

Not only this, but a variety of counterexamples to the corona
problem in SCV have been produced in [FS, S1]; for failure of the
corona problem on a Riemann surface, see [G] for a discussion on
Cole's example. One of the main purposes of this paper is to
provide favorable results on the corona problem when $G$ is a
certain pseudoconvex domain (including the ball) or else a
polydomain (including the polydisc) in $\cc^n$ (presently, to keep
our results optimal, we consider certain domains in $\cc^n$
(including basic domains such as the ball and the polydisc) only,
keeping ourselves away from the complex manifolds (of finite
dimension)).

In fact, after the Carleson's solution, most of the subsequent
attempts to extend it for more general domains in $\cc^n$ were
aimed at solving the B\'{e}zout's equation by various methods
(mostly, by solving the $\bar{\partial}-$ equation); see, for
instance, [KL1-2, K5, T, V]. These approaches required hard
classical analysis methods and has met serious difficulties in the
case when $n > 1$, not leading to the solution to the original
problem; see, for instance, [KL1-2, V]. For example, Varopoulos
realized that his result {\it does not imply} that the corona
problem fails in several variables -- only that the
$\bar{\partial}-$technique with that {\it particular definition of
Carleson measure} fails. In fact, he noted between Thms. 3.1.1 and
3.1.2 why the Carleson's classical approach of solving
$\bar{\partial}_b-$problem fails in the SCV case when one wanted
to prove the corona theorem (see Thm. on p. 272 of [V] wherein the
solution functions are BMO). However, in this paper, we show that
this is {\bf certainly} the case. Indeed, we first generalize the
result of Clos (see Thm. 2.3 below) on the Gleason solvability by
assuming that $G$ is a {\it certain} pseudoconvex domain in
$\cc^n$ (i.e., not just a pseudoconvex domain for which the
$\bar{\partial}-$problem is solvable in $L^\infty$, namely,
$L^\infty-$pseudoconvex domain), and then, using this result, we
establish the corona theorem in Thm. 2.4 below. In fact, in \S 4,
we solve this problem for other familiar Banach spaces of bounded
holomorphic functions on some very familiar domains (including
balls and polydiscs) in $\cc^n$; we then need some well-known
classical functional analytic methods such as the Banach algebra
method and the classical dead-end space method. In \S 2, we first
show an {\it equivalency} of this problem with another significant
problem in the theory of SCV, namely, the Gleason's problem
(1964), which we now describe in fuller details as follows.

\noindent {\bf Gleason's problem.} Let $A(G)$ be the Banach
algebra of all holomorphic functions on a domain $G$ in $\cc^n$
that extend by continuity to the closure of $G$, equipped with the
usual compact open topology given by the sup-norm. Thus, it is a
closed subalgebra of $H^\infty(G)$, and so, for our convenience,
we may also consider $A(G)$ as the Banach algebra of all uniformly
continuous, holomorphic functions on $G$ [A]. The Gleason's
problem was to decide whether the ideal in $A$, where $A$ is
either $H^\infty(G)$ or $A(G)$, consisting of functions vanishing
at a fix point
$\lambda\,=\,(\lambda_1,\,\lambda_2,\,\dots,\,\lambda_n)\,\in\,G$,
is algebraically finitely generated by the shifted coordinate
functions $z_j - \lambda_j$ for $j\,=\,1,\,2,\,\dots,\,n$ (see
[BF] and other references therein for more details). This is
actually a variant of the original Gleason's problem, posed for
$A(G)$, where $G$ is the unit ball in $\cc^2$ and $\lambda = (0,
0)$ [Gl]. The positive solution in the case of strictly
pseudoconvex domains with smooth boundaries was established in
[H2]. It has been even strengthened in [J] as follows: The
continuous linear operators
$$T_j : A(G) \rightarrow
A(G),\;f\,\mapsto\,T_j(f),\;j\,=\,1,\,2,\,\dots,\,n$$ solving the
Gleason's problem exist, which means that $$f(z) -
f(\lambda)\,=\,\sum_{j = 1}^n(z_j -
\lambda_j)(T_j(f))(z),\;z\,\in\,G, f\,\in\,A(G).$$ For polydomains
analogous decompositions are easily available by fixing all but
one variable method.

This implies, in particular, that the part of the spectrum
$M(A(G))$ in the fibers over $G$ is mapped in a bijective way to
$G$ by the Gel'fand transform $(\hat{z}^1,\,\dots,\,\hat{z}^n)$ of
the coordinate functions (see also [A]). By the result of [HS],
the entire spectrum $M(A(G))$ is mapped bijectively onto
$\overline{G}$ where $G$ is a smooth pseudoconvex domain in
$\cc^n$ (it also suffices to consider a pseudoconvex domain $G$
with a Stein neighbourhood basis; see [Ro]). Hence the Gel'fand
topology coincides with the Euclidean topology on $\overline{G}$
and $G$ is open in $M(A(G))$.

\noindent {\bf Necessary and sufficient background.} This paper
contains a continuation of the work begun in [P4], and we shall
feel free to use the terminology and conventions established
there. However, this current work is specifically concerned with
an equivalency of the corona problem and Gleason's problem in the
theory of SCV. As is seen in [P4], we give sufficient conditions
for the existence of (local) analytic structure in the spectrum of
a commutative Fr\'{e}chet algebra by studying the ideal structure
of the algebra. This extends the result of Gleason to finitely
generated ideals in Fr\'{e}chet algebras, answering affirmatively
a question posed in [Ca] for Fr\'{e}chet algebras. As a
consequence, locally Stein algebras are completely characterized
by intrinsic properties within the class of Fr\'{e}chet algebras,
and, as an application of this characterization, an affirmative
answer to the Gleason's problem for such algebras is provided
recapturing all the classical results on the Gleason's problem in
the theory of SCV. In fact, we remark that complex analysts may
also find these sufficient conditions interesting from
applications point of view; for example, see Thms. 4.2, 4.3 and
Cor. 4.4 of [P4], and references (and their reviews, too) to [Gl]
in MathSciNet (MR0159241 (28 \#2458)). As we will see, this is,
indeed, the case. In fact, as is well-known, in the context of
Banach algebras $A(G)$ and $H^\infty(G)$ (and more generally,
Banach algebras of bounded holomorphic functions on $G$, where $G$
is a domain in $\cc^n$), the knowledge of the fiber structure is
useful in the study of Gleason parts which was motivated by the
search for analytic structure in the spectrum.

We recall that Banach algebras satisfying Lem. 3.1 of [P4] (in
particular, the polydisc algebra $A(D^n)$, the ball algebra
$A(B^n)$ and the algebra $H^\infty(G)$) are locally Stein-Banach
algebras. In \S 3, we show that for $\alpha \,\in\,(0, 1]$, Banach
subalgebras $\textrm{Lip}_H^{\alpha}(G)$ of $A(G)$ (and so, of
$H^\infty(G)$ as well) whose Gleason solvability we want to
establish, are also locally Stein-Banach algebras. Further, we
recall that a domain $G$ has the Gleason $A$-property if the
Gleason's problem has an affirmative solution at all points of $G$
for a locally Stein algebra $A$. Clos calls it the Gleason
solvability property of a domain $G$ for $H^\infty(G)$. In fact,
the Gleason solvability is sometimes referred to as Hefer's
condition (as developed by Hefer in 1940). For more information on
the Hefer's condition, see [K3, Ch. 5]. Many of the methods used
in the study of Hefer's condition involve the integral
representation techniques developed in [H1, R]. The Gleason
solvability was first proved by Leibenson, albeit informally [Ru].
The Gleason solvability for holomorphic Bergman spaces,
holomorphic mixed-norm spaces, and the holomorphic Bloch spaces
were studied in [Ru, L, O, RS1, Z]. The domains considered were
initially the unit ball, but then were generalized to strongly
pseudoconvex domains with $C^2$-smooth boundary. There is also
considerable interest on other function spaces and domains,
including the weighted Bergman spaces on egg shaped domains, as
seen in [RS2]. See [Hu] for some results on the Gleason
solvability on the harmonic mixed-norm and Bloch spaces of bounded
strongly pseudoconvex domains with smooth boundary.

Moreover, we recall that a subspace $X$ of $M(A)$ has the Gleason
$A$-property if the problem has an affirmative solution at all
points of $X$. In Cor. 4.4 of [P4], we show that the dense open
subspace $Y$ has the Gleason $A$-property, where $A$ is a
semi-simple locally Stein algebra; in particular, $A$ can be
either $H^\infty(G)$ or $A(G)$, where $G$ is a domain in $\cc^n$
containing the origin. Thus we have given an abstract touch to the
affirmative solution of the Gleason's problem, and so this
abstract method recaptures the classical results obtained by
Beatrous Jr. [Be], Forn\ae ss and \O vrelid [FO], Kerzman and
Nagel [KN], Lieb [Li], Noell [N], and Backlund and
F\"{a}llstr\"{o}m [BF] (see other references therein for a list of
papers on the Gleason's problem), that is, a domain $G$ in $\cc^n$
has the Gleason $A$-property, where $A$ is either $H^\infty(G)$ or
$A(G)$ and $G$ is a (strictly or weakly) pseudoconvex domain in
$\cc^n$ with various boundary conditions.

Concerning the corona problem, as far as we know, almost all
attempts in the past were either relying on the solution of the
B\'{e}zout equation or relying on the use of the Koszul complex
(an idea developed for the computation of the cohomology of a Lie
algebra [K5, T]) while applying the $\bar{\partial}-$method. As we
shall see, we find it necessary (mathematically, ``sufficient") to
consider the Gleason solvability while solving the corona problem.
On the other hand, although we have not spelled it out explicitly,
but have considered the corona theorem as a hypothesis while
solving the Gleason's problem (see Thm. 4.3 and Cor. 4.4 of [P4]).
Moreover, Kerzman and Nagel used the Koszul complex while solving
the Gleason's problem for certain function algebras on certain
domains in $\cc^n$. Now, all this developments stem to think on an
equivalency of these two problems, as well as to give a short
elegant proof, by using the (maximal) ideal structure of the
Banach algebras $A(G)$ and $H^\infty(G)$ alone (that means we try
to keep ourselves away from the analytical aspects such as Gleason
parts of the spectra of these algebras).

One of our main goals in this paper is to extend Thm. 4 of [Cl] in
\S 2 as follows. First, we show that if $G$ is a pseudoconvex
domain in $\cc^n$, then it has the Gleason $A(G)$-property in
Prop. 2.2. Using this, we show that $G$ has, indeed, the Gleason
$H^\infty(G)$-property (which, further, shows that $G$ is open in
$M(H^\infty(G))$ w.r.t. Gel'fand topology). Then, we bring in
Arens' result (in the Banach case) to show that the B\'{e}zout's
equation is solvable, establishing the corona problem for such
domains in $\cc^n$, including balls and polydiscs.

In [CI], it was not evident what relationship there was between
the condition that $G$ be open in $M(H^\infty(G))$ (or $M(A(G))$)
w.r.t. Gel'fand topology and the condition that $G$ be an
$L^\infty$-pseudoconvex domain. However, Clos showed the one-way
implication (see Thm. 4 of [Cl] and the fact, discussed in [CI],
that $G$ is open in $M(H^\infty(G))$ if it has the Gleason
$H^\infty(G)$-property). One reason for interest in this
discussion stems from an application to the commuting Toeplitz
operators problem on Bergman spaces of pseudoconvex domains in
$\cc^n$ (see [Cl, p. 2] for his comment). This issue will be
addressed in a future paper.

It is worthwhile mentioning some work that is somewhat
(tangentially) related to the work given here (see [ACG, ADLM,
AFGM, CGJ, Re] for details). In fact, their work is specifically
concerned with the existence of (local) analytic structure in
spectra of certain algebras of bounded holomorphic functions on
the open unit ball $B_X$ of a complex Banach space $X$; e.g., they
study the Gleason parts of spectra of these algebras when $X \,=\,
c_0$ in [ADLM]. To establish a relation between their work with
the current work, we raise two questions in 4.4 of [P4].
\section{Equivalency of the two problems.}
As discussed above, we first show that Thm. 4 of [Cl] holds true
for the Banach algebra $A(G)$, where $G$ is a pseudoconvex domain
in $\cc^n$. This is possible because of Thm. 5 of [Cl] and the
following
\begin{lem} \label{A(G)_subset} Let $G$ be a pseudoconvex domain
in $\cc^n$. Then $A(G)$ is a Banach subspace of the Bergman space
$A^2(G)$.
\end{lem}
\noindent {\bf Proof.} Since $A(G)$ is the Banach algebra of
uniformly continuous, holomorphic functions on a pseudoconvex
domain $G$ in $\cc^n$, we have that $A(G) \subset H^\infty(G)
\subset \textrm{Hol}(G)$ w.r.t. the usual compact-open topology.
This implies that
$$A(G) \bigcap L^2(G)\subset H^\infty(G) \bigcap L^2(G) \subset
\textrm{Hol}(G)\bigcap L^2(G) = A^2(G),$$ where $A^2(G)$ is the
Bergman space of a pseudoconvex domain $G$. The polynominals in
$z_1,\,\dots,\,z_n$ are dense in $A^2(G)$ since $A^2(G)$ is
isometrically isomorphic to the weighted space
$\ell^\infty((\zz^+)^{n}, \frac{1}{|N| + 1})$ (see [P3] for the
analogous notion of the Beurling-Banach algebra
$\ell^1((\zz^+)^{n}, \omega)$ of weight type). Since each function
of $A(G)$ is uniformly continuous on $G$, it extends uniquely and
in a norm preserving way to a continuous function on
$\overline{G}$ (see \S 1), it is obviously square-integrable on
$\overline{G}$. Hence $f\,\in\,L^2(G)$ whenever $f\,\in\,A(G)$,
and so, $f\,\in\,A^2(G)$. Thus $A(G)$ is a Banach subspace of the
Bergman space $A^2(G)$ (certainly, $A(G)$ is not a closed subspace
of $A^2(G)$ for an obvious reason).$\hfill \Box$

Thus we have the following proposition, extending results of [Cl,
KN].
\begin{prop} \label{Thm 4_A(G)} Let $G$ be a pseudoconvex domain
in $\cc^n$. Suppose $f \in A(G)$ and $f(\a) = 0$ for some $\a \in
G$. Then there exists $f_1,\,f_2,\,\dots,\,f_n$ in $A(G)$ so that
$f\,=\,\sum_{j = 1}^n(z_j - \alpha_j)f_j$. In particular, $G$ is
Gleason solvable for $A(G)$.
\end{prop}
\noindent {\bf Proof.} Since $f\,\in\,A(G)$, $f\,\in\,A^2(G)$ by
Lem. \ref{A(G)_subset}. Then, by [Cl, Thm. 5], there exists
$f_1,\,f_2,\,\dots,\,f_n$ in $A^2(G)$ so that $f\,=\,\sum_{j =
1}^n(z_j - \alpha_j)f_j$. Since $f\,\in\,A(G)$ and shifted
coordinate functions $z_j - \alpha_j$ for $j\,=\,1, \,2,
\,\dots,\,n$, are in $A(G)$, it is clear that
$f_1,\,f_2,\,\dots,\,f_n$ are also in $A(G)$. $\hfill \Box$

As an application of the above result, we extend Thm. 4 of [Cl] as
follows. Below, by a ``pseudoconvex" domain in $\cc^n$, we mean a
smooth (strictly) pseudoconvex domain in $\cc^n$ [HS] (it also
suffices to consider a pseudoconvex domain with a (strong) Stein
neighbourhood basis [Ro, Sa]; e.g., Stein domains in $\cc^n$).
Moreover, we also remark that if $G$ is a relatively compact
domain in a Riemann surface $M$, then the maximal ideal space
$M(A(G))$ is $\overline{G}$, and so is true if $G$ is relatively
compact domain such that $\overline{G}$ is an $S_\delta$ set
(equivalently, $G$ has a Stein neighbourhood basis) in a complex
analytic manifold $M$ (see [Ro, Thm. 2.12]).
\begin{thm} \label{Thm 5_A(G)} Let $G$ be either a ``pseudoconvex"
domain or else a polydomain in $\cc^n$. Then $G$ is
Gleason solvable for $H^\infty(G)$.
\end{thm}
\noindent {\bf Proof.} First, we recall that the Banach algebra
$A(G)$ is a closed subalgebra of the Banach algebra $H^\infty(G)$
w.r.t. the compact-open topology induced by the sup-norm such that
the polynomials in $z_1,\,\dots,\,z_n$ are dense in $A(G)$, but
not in $H^\infty(G)$. In fact,
$\textrm{in}\,:\,A(G)\,\hookrightarrow\,H^\infty(G)$ is a
continuous injective homomorphism. So, the adjoint spectral
mapping $\textrm{in}^*\,:\,M(H^\infty(G))\,\rightarrow
\,M(A(G))\,=\,\overline{G}$ is surjective and continuous w.r.t.
the Gel'fand topology. Indeed, it's restriction is one-to-one over
$G$, and thus, the inverse mapping maps $G$ homeomorphically onto
an open subset of $M(H^\infty(G))$ (see [A, ACG, CGJ]). In fact,
we remark that, in [ACG], the one dimensional case of $G$ being
the open unit disc, was explained [Ho], and claimed that the same
facts are probably true and can be proved at least whenever $G$ is
the open unit ball or the polydisc in $\cc^n$; for proof of the
latter claims, see [A], which may be applied to the case when $G$
is either a pseudoconvex domain or a polydomain in $\cc^n$.

Next, we show that the maximal ideal $M_\alpha$, $\alpha\,\in\,G$,
of $H^\infty(G)$ is, indeed, algebraically finitely generated. For
this, we note that the maximal ideal $M_\alpha$, $\alpha\,\in\,G$,
of $A(G)$ is algebraically finitely generated by Prop.~\ref{Thm
4_A(G)}. Recall that when $G$ is a polydomain in $\cc^n$, the
Gleason solvability is easily available through the analogous
decomposition, discussed in \S 1, by fixing all but one variable
method. Then, by the discussion above, it's inverse image, denoted
by $M_\alpha$ again, is also a maximal ideal of $H^\infty(G)$. To
show that it is, indeed, finitely generated, we remark that the
finitely generated maximal ideal $M_\alpha$ of $A(G)$ is, in fact,
an ideal of $H^\infty(G)$. By the Zorn's lemma, the maximal
element of the family $\mathcal I$ containing finitely generated
ideals of $H^\infty(G)$ consisting of functions in $H^\infty(G)$
vanishing at $\alpha$ in $G$ (note that $\mathcal I$ is non-empty,
since $M_\alpha\,\in\,\mathcal I$), is an algebraically finitely
generated maximal ideal. In fact, this maximal element is not only
a primary ideal, but is obviously contained in the maximal ideal
$M_\alpha$ of $H^\infty(G)$. From these facts, we conclude that
the maximal element is, indeed, the maximal ideal $M_\alpha$ of
$H^\infty(G)$. $\hfill \Box$

We remark that there are smooth pseudoconvex domains in $\cc^n$
for which $L^\infty-$estimates for $\bar{\partial}$ do not hold
(see [B, FS, S2]), and thus, the above theorem is a significant
generalization of Thm. 4 of [Cl]. In order to establish an
equivalency between the corona problem and Gleason's problem, we
first show that the Gleason solvability, obtained above, implies
an affirmative answer to the corona problem in the following
theorem. We remark that there are smooth pseudoconvex domains in
$\cc^n, \;n\,=\,2,\,3$, for which the corona theorem fails [FS,
S1]. However, below, by a {\bf pseudoconvex} domain in $\cc^n$, we
mean a relatively compact, smooth (strictly) pseudoconvex domain
in $\cc^n$ such that its closure is an $S_\delta$ set (i.e., it
has a (strong) Stein neighbourhood basis) [Ro, \S 5, Defn. 2] and
[Sa, Introduction]. Such a pseudoconvex domain in $\cc^n$
possesses (strictly) pseudoconvex boundary and is of type
$HL^\infty$; i.e., it is a $H^\infty$ domain of holomorphy [K6,
Prop. 6]. Moreover, we also remark that if $G$ is a relatively
compact, smooth (strictly) pseudoconvex domain in a complex
manifold $M$ such that $\overline{G}$ is an $S_\delta$ set
(equivalently, $G$ has a Stein neighbourhood basis), then it
possesses a (strictly) pseudoconvex boundary.
\begin{thm} \label{Corona} Let $G$ be either a {\bf pseudoconvex} domain or else a
polydomain in $\cc^n$. Then the corona problem holds true in
$H^\infty(G)$. That is, the spectrum $M(H^\infty(G))$ is the
Gel'fand closure of (the canonical image of) $G$.
\end{thm}
\noindent {\bf Proof.} By Thm.~\ref{Thm 5_A(G)}, $G$ is Gleason
solvable. To show that the corona problem has an affirmative
solution for $G$, we use its equivalent formulation, discussed in
\S 1, namely, a finite {\it corona data} implies the existence of
the solution of the B\'{e}zout equation. In fact, we show that if
there does not exist a solution of the B\'{e}zout equation, then
there are no finite corona data available. To this effect, we
recall that there is a one-to-one correspondence between the
maximal ideals and the elements of the spectrum of the Banach
algebra $A$. As an easy consequence of this fact, we further know
that the algebraically finitely generated ideal, generated by
$f_1, \dots, f_n \in A$, is proper if and only if
\begin{equation}\label{ArensBanach} \{\phi\,\in\,M(A)\,:\,\phi(f_j)\,=\,0,\;j\,=\,1,\,\dots,\,n\}\,\neq\,\emptyset.
\end{equation}
In our case, since $G$ is Gleason solvable, it turns out that
$M_\alpha$, $\alpha\,\in\,G$, is algebraically finitely generated
maximal ideal if and only if the equation~\ref{ArensBanach} holds
if and only if there are no $g_1,\,\dots,\,g_n\,\in\,H^\infty(G)$
such that $\sum_{j = 1}^nf_jg_j\,=\,1$, where $f_j,\; j
\,=\,1,\,\dots,\,n,$ are the generators of $M_\alpha$. Thus there
does not exist a solution of the B\'{e}zout equation if and only
if there are NO finite corona data available (that means,
$\alpha\,\in\,G$ is such that for every $\delta\,>\,0$, $\sum_{j =
1}^n|f_j(\alpha)|\,<\,\delta$, which implies that $\sum_{j =
1}^n|f_j(\alpha)|\,=\,0$ if and only if $f_j(\alpha)\,=\,0$).
Hence the corona problem holds true in $H^\infty(G)$, where $G$ is
either a {\bf pseudoconvex} domain or a polydomain in $\cc^n$.
$\hfill \Box$

\noindent {\bf Remarks A.} Here are some remarks on our
Banach-algebraic technique in the above proof.

\noindent 1. In connection with Krantz's Thm. 1 and Rem. 2 [K5],
we note that if we apply the contrapositive method to prove Thm. 1
of Krantz, then we must show that failure of the corona theorem
would imply that either there are NO corona data or, the zeros of
the corona data are NOT separated. However, it is clear that there
are NO corona data if and only if the zeros of the functions in
corona data are NOT separated, as explained in the proof above.
That means, one may drop the geometric hypothesis, imposed in Thm.
1 of Krantz, in order to further simplify that theorem. Now, if we
apply the contrapositive method to prove the simplified Thm. 1 of
Krantz, then we have the corona theorem (and that is what we have
done in the above proof for a finite corona data).

\noindent 2. Our Banach-algebraic technique in the above proof
also works for an infinite corona data (i.e., with infinitely many
generators $f\,=\,(f_1,\,f_2,\,\dots)$), taken appropriately (see
[DKSTW, \S 3.5] for the disc in $\cc$) due to [Go, Brook's Thm.,
p. 136]; note that we not only recover Thm. 5 of Rosenblum and
Tolokonnikov in some sense, but extend it for more general $G$ in
$\cc^n$.

\noindent 3. It is interesting to note that author very lately
(i.e., while finalizing the initial draft) came to know about the
{\bf ideal problem}, explained in the history of the corona
problem (see [DKSTW, \S 7]). So, our Banach-algebraic technique is
purely based on routine results from the Banach algebra theory as
well as author's recent work on the Gleason's problem in [P4].
However, some more comments are in order. It is easy to see that
our method here is similar to that of Xiao's theorem on the
multipliers of the Dirichlet space $\cal D$ in the unit disc (see
[DKSTW, p. 22] or, [X]) where it is discussed that the corona data
is necessary for the corona theorem to hold (see also a remark on
an alternative proof which also deals with infinitely many
generators; p. 23 of [DKSTW]), and in Thm. 4.1 below, we show that
the corona data alone is sufficient for the corona theorem to hold
(see [KL1, K5]). Further, it is easy to see that we have, in fact,
solved the Rubel and Shields' question for the function $\phi
\,\equiv\,1$ in $H^\infty(G)$, where $G$ is as above (including
balls and polydiscs in the below corollary) [DKSTW, \S 7,
Questions 1, 2 and 3]; for Questions 2 and 3, we take
$h\,\equiv\,1$ on $[0, \,\infty]$ and $\phi\,\equiv\,1$ in the
above theorem (note that, for Question 3 to answer, we consider a
maximal ideal $M_\a,\;\a\,\in\,G$, generated by
$f_1,\,\dots,\,f_n$ in the above theorem).

As a result, since the ball is (essentially the only) basic
example of a pseudoconvex domain in $\cc^n$ of our kind, we have
the following
\begin{cor} \label{Corona_Special} Let $G$ be either an open unit ball or a polydisc in
$\cc^n$. Then the corona problem holds true in $H^\infty(G)$.
$\hfill \Box$
\end{cor}
We have already shown in Thm.~\ref{Thm 5_A(G)} that if $G$ is
either a smooth (strictly) pseudoconvex domain in $\cc^n$ or a
pseudoconvex domain in $\cc^n$ with a (strong) Stein neighbourhood
basis (e.g., Stein domains in $\cc^n$) or a polydomain in $\cc^n$,
then $G$ is Gleason solvable. However, as a part of our goal to
establish an equivalency between the two problems, in the converse
direction, we have the following
\begin{thm} \label{Gleason} Let $G$ be either a
relatively compact, smooth (strictly) pseudoconvex domain in
$\cc^n$ with a (strong) Stein neighbourhood basis or else a
polydomain in $\cc^n$. Then $G$ is Gleason solvable provided that
the corona problem holds true in $H^\infty(G)$.
\end{thm}
\noindent {\bf Proof.} Let the corona problem hold true in
$H^\infty(G)$. It is easy to show that $H^\infty(G)$ is a
semisimple locally Stein-Banach algebra by applying proof of (ii)
implied (i) of Thm. 4.1 of [P4]. Hence, by Thm. 4.3 of [P4] (when
applied to $H^\infty(G)$), we see that the Gleason's problem is
solvable, because $G$ is dense in the spectrum $M(H^\infty(G))$.
$\hfill \Box$

\noindent{\bf Remarks B.} 1. By Thms. 2.4 and 2.6, we see that the
two significant problems from the theory of SCV are equivalent. In
fact, the corona problem has a topological nature (i.e., whether
the domain $G$ is dense in the spectrum $M(H^\infty(G))$;
equivalently, whether the corona is empty) whereas the Gleason's
problem has an algebraic nature (i.e., whether a maximal ideal
$M_\a,\, \a\,\in\,G$, is algebraically finitely generated in
$H^\infty(G)$; equivalently, whether $G$ has the Gleason
$H^\infty(G)-$property). The corona problem has also an algebraic
nature (i.e., given a finite corona data whether the B\'{e}zout's
equation holds true). In the theory of automatic continuity, we
are interested in the interplay between the topological nature and
the algebraic nature of the underlying structure (see [DPR, P1,
P3]).

\noindent 2. Since we use the corona theorem as a hypothesis in
the above theorem, we obtain the Gleason solvability for a limited
class of domains in $\cc^n$ (see the preceding paragraph to Thm.
2.4) than we obtain the same for a larger class of domains in
$\cc^n$ in Thm. 2.3.
\section{Gleason problem for Banach spaces of holomorphic functions.}
{\bf 3.1} [FO] is a good reference; in the final paragraph of \S
1, the authors state the main theorem can still be proved by
replacing $A(G)$ by various holomorphic H\"{o}lder- and
Lipschitz-spaces and by replacing the coordinate functions as the
generators by arbitrary generators of the maximal ideal in these
spaces. They mention a future paper, but, as far as we know, it
has never been published. In this connection, we can consider the
below application of Thm. 4.3 of [P4], which establishes the
similar claim, but for semi-simple locally Stein algebras. Indeed,
we show that for the exponent $\alpha\,\in\,(0, 1]$, the Banach
algebras $\textrm{Lip}_H^\alpha(G, d)$ of all bounded,
complex-valued, $\alpha-$Lipschitz continuous (aka H\"{o}lder
continuous), holomorphic functions on $(G, d)$, where $G$ is a
bounded domain in $\cc^n$ and $d$ is the usual metric on $\cc^n$,
is a semisimple locally Stein-Banach algebra under the usual
``$\alpha-$Lipschitz" norm (defined below). Thus, if $\alpha = 1$,
then it is, indeed, the Banach algebra $\textrm{Lip}_H(G, d)$ of
all bounded, complex-valued, Lipschitz continuous, holomorphic
functions on $(G,\, d)$. Below, we establish the corona theorem in
a stronger sense, and as an application of this theorem, we
establish the Gleason solvability of certain bounded pseudoconvex
domains in $\cc^n$ for these Banach algebras (see Lem. 3.2 and
Thm. 3.13), including an affirmative solution to the question left
open by Forn\ae ss and \O vrelid.

Moreover, for $\a\,=\,1$, the Banach algebra
$\textrm{Lip}(X,\,d^{'})$ of all Lipschitz continuous functions on
the metric space $(X, \,d^{'})$ was studied by Sherbert in
[Sh1-2], and for each $\a\,\in\,(0,\,1)$, the Banach algebra
$\textrm{lip}^\a(X,\,d^{'})$ was studied in [Sh2]. As far as we
know, for each $\a\,\in\,(0,\,1]$, the Banach algebra
$\textrm{Lip}_H^\a(G,\,d)$ was never studied, which is, indeed, an
important algebra compared to its superalgebra
$\textrm{Lip}^\a(G,\,d)$, especially when it comes to further
development of the theory of Banach algebras of holomorphic
functions (remark that for each $\a\,\in\,(0,\,1)$, the Banach
algebra $\textrm{Lip}^\a(X,\,d^{'})$ of all Lipschitz continuous
functions of order $\a$ on the metric space $(X,\,d^{'})$ was also
{\bf never studied} by Sherbert). So, in this subsection, we study
the Banach algebras $\textrm{Lip}_H^\a(G,\,d)$ in view of the two
problems in the theory of SCV (e.g., Gleason's problem, discussed
in [P4]). Thus our results extend results from [Sh1-2]. We may
refer to these algebras as Lipschitz algebras of holomorphic
functions. We examine some basic properties of these Banach
algebras.

First, we will discuss the Gel'fand theory of
$\textrm{Lip}_H^\a(G,\,d)$, where $\a\,\in\,(0,\,1]$, and show
that every commutative semisimple Banach algebra $A$ of
holomorphic functions (i.e., for every $x\,\in\,A$, $\hat{x}$ is
holomorphic on the interior of its maximal ideal space $M(A)$
w.r.t. the Gel'fand topology), is isomorphic to a subalgebra of
the Lipschitz algebra $\textrm{Lip}^\a(M(A),\,\sigma)$ of order
$\a$, where $\sigma$ is the metric on the maximal ideal space
$M(A)$ inherits from as a subset of the dual space $A^*$ as
defined below. This representation is drawn from the Gel'fand
representation. We use the metric topology in place of the usual
Gel'fand topology of $M(A)$. We also show that this isomorphism is
onto if and only if $A \,=\,\textrm{Lip}_H^\a(G,\,d)$ for a
bounded domain $G$ in $\cc^n$. That is, the image of
$\textrm{Lip}_H^\a(G, \,d)$ under the Gel'fand mapping results in
strictly those functions in $\textrm{Lip}^\a(M(A), \,\sigma)$ that
are continuous on $M(A)$ and holomorphic on $\textrm{Int}\,M(A)$
w.r.t. the Gel'fand topology (and so, these functions are also
continuous on $M(A)$ and holomorphic on $\textrm{Int}\,M(A)$
w.r.t. the $\sigma-$topology). From this result, it is clear that
these algebras are semisimple, and so, they have a unique complete
norm.

We recall that we used the {\it corona theorem} as a hypothesis to
solve the Gleason's problem for locally Stein algebras. However,
it was not evident what relationship there is between the {\it
corona theorem} and the {\it solvability of the Gleason's
problem}, i.e., it was not obvious whether either of these
problems implies the other. This issue of equivalency of both the
problems is addressed in \S 2 above, solving the most prestigious,
long-standing, original corona problem for the Banach algebra
$H^\infty(G)$, where $G$ is a certain bounded domain in $\cc^n$
(including the unit ball and the polydisc). The solvability of
both the problems for various Banach spaces of bounded holomorphic
functions (e.g., holomorphic mean Besov-Lipschitz spaces,
Besov-Lipschitz spaces, Hardy-Sobolev spaces and a weighted
Bergman space) on various types of bounded domains in $\cc^n$
(including the unit ball and the polydisc) is discussed in the
below (sub)sections. Finally, we show two significant differences
between these two classes of Lipschitz algebras: for each
$\a\,\in\,(0,\,1]$, $\textrm{Lip}^\a(X,\,d^{'})$ is a regular,
self-adjoint (closed under complex conjugation) Banach algebra
whereas $\textrm{Lip}_H^\a(G,\,d)$ is neither regular nor
self-adjoint.\\

\noindent {\bf Gel'fand theory of Lipschitz algebras of
holomorphic functions.} For each $\alpha \in (0, 1]$, let
$\textrm{Lip}_H^\alpha(G, d)$ denote the collection of all
bounded, complex-valued, $\alpha-$Lipschitz continuous (aka
H\"{o}lder continuous), holomorphic functions defined on the
bounded domain $(G, d)$ in $\cc^n$, where $d$ is the usual metric
on $\cc^n$ induced by either the max norm or the square norm.
Thus, for a fixed $\a\,\in\,(0,\,1]$, $\textrm{Lip}_H^\a(G,\, d)$
consists of all $f$ defined on $(G,\, d)$ such that both
$$\|f\|_\infty\,:=\,\sup\{|f(z)|\,:\,z\,\in\,G\}$$ and $$\|f\|_\alpha\,:=\,\sup\{\frac{|f(z) - f(w)|}{\|z - w\|^\alpha}\,:\,z, w\,\in\,G, z\,\neq
\,w\}$$ are finite. It is easy to see that, with the norm
$\|\cdot\|_{\textrm{Lip}^\alpha}$, defined by
$\|f\|_{\textrm{Lip}^\alpha}\,=\,\|f\|_\infty\, +\, \|f\|_\alpha$,
$\textrm{Lip}_H^\alpha(G,\, d)$ is a commutative Banach algebra
with identity (note that another equivalent algebra norms are also
possible, e.g., if we take maximum of $\|f\|_\infty$ and
$\|f\|_\a$ or, if we take supremum of $\{|f(z_0)|,\,\frac{|f(z) -
f(w)|}{\|z - w\|^\a}\,:\,z_0,\,z,\,w\,\in\,G, \,z\,\neq\,w\}$;
also, we will see that there are several equivalent norms possible
on the Banach space $\textrm{Lip}_H^\alpha(G,\,
d)\,=\,H\Lambda(\alpha/n,\,\infty,\,\infty)\,=\,\Lambda_\alpha\,=\,\Lambda_\alpha^{\infty,
\,\infty}$ due to various characterizations discussed in [F, Ch,
CZ, JP]). In fact, it is a closed subalgebra of a Banach algebra
$\textrm{Lip}^\a(G,\,d)$ of $\alpha-$Lipschitz (aka H\"{o}lder)
continuous functions (resp., Banach algebra $\textrm{Lip}(G,\,d)$
of Lipschtiz continuous functions, when $\a\,=\, 1$). We may
sometimes refer to such an algebra as a H\"{o}lder algebra of
bounded holomorphic functions (resp., Lipschitz algebra of bounded
holomorphic functions, when $\a\, =\, 1$). For the exponent
$\alpha\,\in\,(0, 1]$, it is evident that
$$\textrm{Lip}_H(G)\,\subset
\,\textrm{Lip}_H^\alpha(G)\,\subset\,A(G)\,\subset\,H^\infty(G),$$
where $G$ is a bounded domain in $\cc^n$.

It will be assumed throughout the subsection that $(G,\, d)$ is a
relatively compact metric space in $\cc^n$ unless otherwise
specified. Hence we see that since each element of
$\textrm{Lip}_H^\alpha(G,\, d)$ is strongly uniformly continuous
on $(G,\, d)$ (resp., uniformly continuous on $(G,\, d)$, when
$\a\, =\, 1$), it extends uniquely, continuously and in norm
preserving way to an element of $\textrm{Lip}^\a(\overline{G},\,
d)$. Thus, as Banach algebras, $\textrm{Lip}_H^\alpha(G,\, d)$ and
its image in $\textrm{Lip}^\alpha(\overline{G},\, d)$ are
isometrically isomorphic (equivalently,
$\textrm{Lip}_H^\alpha(G,\,d)$ may also be considered as a closed
subalgebra of a Banach algebra $\textrm{Lip}^\alpha(\overline{G},
\,d)$).

Let $A$ be a commutative semisimple Banach algebra of holomorphic
functions with norm $\|\cdot\|_A$. By the standard Banach space
theory, the maximal ideal space $M(A)$ lies on the unit sphere of
the dual space $A^*$. Thus, as a subset of $A^*$, $M(A)$ inherits
the Gel'fand topology, as well as the relative norm or metric
topology. The Gel'fand theory of commutative Banach algebras uses
the Gel'fand topology. When $A$ has an identity, $M(A)$ is compact
w.r.t. the Gel'fand topology. For each $f\,\in\,A$, the function
$\hat{f}$ is defined on $M(A)$ by
$\hat{f}(\phi)\,=\,\phi(f),\;\phi\,\in\,M(A)$. Each $\hat{f}$ is
continuous on $M(A)$ and holomorphic on $\textrm{Int}\,M(A)$
w.r.t. the Gel'fand topology. In fact, the Gel'fand topology is
the weakest topology on $M(A)$ such that the family
$\{\hat{f}\,:\,f\,\in\,A\}$ is continuous on $M(A)$ (and
holomorphic on $\textrm{Int}\,M(A)$). Let $C(M(A))$ will denote
the algebra of complex-valued continuous functions on $M(A)$
w.r.t. the Gel'fand topology, given with the sup norm. Then the
Gel'fand mapping $f\,\rightarrow\,\hat{f}$ is a continuous
isomorphism of $A$ into $C(M(A))$.

We now consider the metric topology of $M(A)$. The metric $\sigma$
on $M(A)$ induced by the norm $\|\cdot\|_A^*$ of the dual space
$A^*$, is defined by
$$\sigma(\phi,\,\psi)\,=\,\|\phi -
\psi\|_A^*,\;\phi,\,\psi\,\in\,M(A).$$ In terms of the functions
$\hat{f},\;f\,\in\,A$, we may express the metric $\sigma$ by
$$\sigma(\phi,\,\psi)\,=\,\sup\{|\hat{f}(\phi) -
\hat{f}(\psi)|\,:\,\,f\,\in\,A,\,\|f\|_A\,\leq\,1\},\;\phi,\,\psi\,\in\,M(A).$$
The metric topology of $M(A)$ is stronger than the Gel'fand
topology. Therefore, since $M(A)$ is closed in $A^*$ w.r.t. the
weak* topology, it is also closed w.r.t. the metric topology.
Hence, $(M(A), \,\sigma)$ is a complete metric space. The metric
$\sigma$ is bounded as $M(A)$ lies on the unit sphere of $A^*$.

With this metric space $(M(A),\, \sigma)$ for a fixed exponent
$\alpha\,\in\,(0, 1]$, we form the Lipschitz algebra
$\textrm{Lip}^\a(M(A), \,\sigma)$ with norm
$\|\cdot\|_{\textrm{Lip}^\a}$ defined by
$$\|g\|_{\textrm{Lip}^\a}\,=\,\|g\|_\infty\,+\,\|g\|_\a,\;g\,\in\,\textrm{Lip}^\a(M(A),\,\sigma).$$
We now show that the Gel'fand mapping takes $A$ into this
Lipschitz algebra to give a ``Lipschitz representation of order
$\a$" in the following
\begin{prop} \label{Sherbert_Lip representation} Let $A$ be a
unital commutative semisimple Banach algebra of holomorphic
functions. Then the Gel'fand mapping is a continuous isomorphism
of $A$ onto a subalgebra of $\textrm{Lip}^{\alpha}(M(A),
\,\sigma)$. Furthermore, for each $f\,\in\,A$,
$\|\hat{f}\|_\a\,\leq\,\|f\|_A$ and
$\|\hat{f}\|_\infty\,\leq\,\|f\|_A$.
\end{prop}
\noindent {\bf Proof.} First, for a fixed exponent
$\a\,\in\,(0,\,1]$, we discuss an important subset
$\textrm{lip}^\a(M(A),\,\sigma)$ of
$\textrm{Lip}^\a(M(A),\,\sigma)$ consisting of all those functions
$f\,\in\,\textrm{Lip}^\a(M(A),\,\sigma)$ with the property that
$$\frac{|f(z) -
f(w)|}{\sigma(z,\,w)^\a}\,\rightarrow\,0\;\textrm{as}\;\sigma(z,\,w)^\a\,\rightarrow\,0.$$
It is certainly a closed subalgebra of $\textrm{Lip}^\a(M(A),
\,\sigma)$; this was established by Mirkil [M] for the case $X\,
=\, [0,\, 2\pi]$ and $d(x,\,y)^\a\, =\,|x - y|^\a,\; 0 < \a < 1$,
but his proof is valid in general.

Next, it is easy to show that for each $\a \,\in\,(0,\, 1)$,
$$\textrm{Lip}(M(A), \,\sigma)\,\subset\,
\textrm{lip}^\a(M(A),\,\sigma)\,\subset\,
\textrm{Lip}^\a(M(A),\,\sigma)$$ and the Banach algebra
$\textrm{lip}^\a(M(A), \,\sigma)$ (and so, $\textrm{Lip}^\a(M(A),
\,\sigma)$ as well) separates the points of $M(A)$ (see [Sh2,
Prop. 1.6]). Further, these inclusions are continuous. Now, by
Prop. 2.1 of [Sh1], the proof is straightforward. Moreover, the
subalgebra $\hat{A}$ of $\textrm{Lip}^{\alpha}(M(A), \,\sigma)$
consists of functions that are continuous on $(M(A),\, \sigma)$
and holomorphic on $(\textrm{Int}\,M(A),\,\sigma)$ because they
are continuous on $M(A)$ and holomorphic on $\textrm{Int}\,M(A)$
w.r.t. the Gel'fand topology. $\hfill \Box$

Next, we study the maximal ideal space $M(A)$ of the Banach
algebra $A\,=\,\textrm{Lip}_H^{\alpha}(G, \,d)$, where
$\a\,\in\,(0,\,1]$, in view of some important properties that are
useful in establishing the corona theorem and the Gleason
solvability for this algebra, as discussed below. We first remark
that it contains the identity element, namely, the function
identically $1$. It also separates the points of $G$ as it also
contains the co-ordinate functions, but it is not self-adjoint. So
each $z \in G$ can be identified with the point evaluation
functional $\phi_z$ in $M(A)$, where $\phi_z(f) = f(z)$. More
precisely, the injection mapping $z \rightarrow \phi_z$ is
one-to-one from $G$ to $M(A)$ and we may regard $G$ as a subset of
$M(A)$. The metric $\sigma$ on $G$ when restricted to $G$, can be
expressed by
$$\sigma(z,\, w)\,=\,\textrm{sup}\{|f(z) -
f(w)|\,:\,f\,\in\,\textrm{Lip}_H^{\alpha}(G,\,
d),\;\|f\|_{\textrm{Lip}^\alpha}\,\leq\,1\},\;z,\,w\,\in\,G.$$

An algebra $B$ of functions defined on a set $X$ is called
inverse-closed if for every function $f\,\in\,B$ satisfying
$|f(x)|\,\geq\,\epsilon\,>\,0$ for all $x\,\in\,X$, then
$f^{-1}\,\in\,B$ as well. Clearly, it is easy to see that for each
$\a\,\in\,(0,\,1]$, $\textrm{Lip}_H^\a(G,\, d)$ is an
inverse-closed algebra (see, e.g., Prop. 1.7 of [Sh2]). We remark
that $\textrm{Lip}_H^{\alpha}(G,\, d)$ is a Banach algebra such
that $\textrm{Lip}_H^{\alpha}(G,\, d)\,\subsetneq\, A(G)$ (and
this inclusion is continuous as $\|f\|_\infty
\,\leq\,\|f\|_{\textrm{Lip}^\alpha})$. Since polynomials in $z_1,
z_2,\dots,z_n$ are dense in $A(G)$ and they are also in
$\textrm{Lip}_H^{\alpha}(G, \,d)$, $\textrm{Lip}_H^{\alpha}(G,
\,d)$ is a dense Banach subalgebra of $A(G)$. Thus,
$\textrm{Lip}_H^{\alpha}(G, \,d)$ is a semisimple, locally
Stein-Banach algebra by Lem. 3.1 of [P4], since the completion of
$\textrm{Lip}_H^{\alpha}(G, \,d)$ under the sup-norm is $A(G)$
while applying proof of (ii) implied (i) of Thm. 4.1 of [P4] in
the Banach case. Moreover, the maximal ideal space of $A(G)$ is
the joint spectrum $\textrm{sp}(z_1,\,\dots,\,z_n)$, and since
$\textrm{Lip}_H^{\alpha}(G,\, d)$ is also inverse-closed in $A(G)$
by [BO, Prop. 1], we have that
$M(A)\,=\,\textrm{sp}(z_1,\,\dots,\,z_n)$. Further, if $G$ is a
smooth (strictly) pseudoconvex domain in $\cc^n$ [HS], or else, if
$G$ is a pseudoconvex domain with a Stein neighbourhood basis [Ro]
(e.g., Stein domains in $\cc^n$), then
$\textrm{sp}(z_1,\,\dots,\,z_n)\,=\,\overline{G}$. Hence $M(A)
\,=\,\overline{G}$ in this case. This gives the following
\begin{lem} \label{BO_spectrum} Let $M(A)$ be the maximal ideal
space of $\textrm{Lip}_H^{\alpha}(G,\, d)$. Then $M(A)
\,=\,\overline{G}$ in the Gel'fand topology provided that either
$G$ is a smooth (strictly) pseudoconvex domain in $\cc^n$, or is a
pseudoconvex domain with a Stein neighbourhood basis. In
particular, if $(G,\, d)$ is relatively compact in $\cc^n$, then
the Gel'fand topology coincides with the $d$-topology of $G$.
$\hfill \Box$
\end{lem}
\noindent {\bf Remarks C.} 1. There are bounded domains in $\cc^n$
such that neither they have a smooth boundary (e.g., Hartogs
triangle) nor they have a Stein neighbourhood basis (e.g., worm
domain) [DF]. For the Hartogs triangle, we may not have $M(A(G)) =
\textrm{sp}(z_1,\,\dots,\,z_n) = \overline{G}$, but for the worm
domain, we have.

\noindent 2. By proving the above lemma, we have, indeed, solved
the corona problem for commutative semisimple Banach algebras
$\textrm{Lip}_H^\a(G,\,d)$, where $\a\,\in\,(0,\,1]$ (here $G$ can
be a relatively compact polydomain in $\cc^n$ as well; see Thm.
3.13 below). In fact, in this case, we have better result in the
sense that the relative Gel'fand topology of $G$ is equivalent to
the $d-$ and $\sigma-$ topologies of $G$ (see Props. 3.3, 3.9 and
Cor. 3.12 below).

\noindent 3. Note that for the equivalency of the Gel'fand
topology with the $d$-topology of $G$, the relative compactness of
$G$ is as such not required. Indeed, we have the following
proposition whose proof we omit (see Prop. 3.3 of [Sh1]).
\begin{prop} \label{prop_3.3} Let $G$ be either a smooth (strictly)
pseudoconvex domain in $\cc^n$, or a pseudoconvex domain with a
Stein neighbourhood basis. Then the relative Gel'fand topology of
$G$ and the $d$-topology of $G$ are equivalent. $\hfill \Box$
\end{prop}
We now turn to the Gel'fand representation of
$A\,=\,\textrm{Lip}_H^{\alpha}(G,\, d)$, where $\a\,\in\,(0,\,1]$.
Prop. 3.1 tells us that the Gel'fand mapping $f\,\mapsto\,\hat{f}$
takes $\textrm{Lip}_H^{\alpha}(G,\, d)$ isomorphically into
$\textrm{Lip}^{\alpha}(M(A),\, \sigma)$, and satisfies
$\||\hat{f}\||_\infty\,\leq\,\|f\|_{\textrm{Lip}^\alpha}$ and
$\||\hat{f}\||_\a\,\leq\,\|f\|_{\textrm{Lip}^\alpha}$ for all
$f\,\in\,\textrm{Lip}_H^{\alpha}(G, \,d)$. These statements pursue
from the general consideration. In the particular case of
$A\,=\,\textrm{Lip}_H^{\alpha}(G,\, d)$, this can be strengthened
in the following
\begin{thm} \label{Sherbert_PropositionCorollary} The Gel'fand
mapping $f\,\mapsto\,\hat{f}$ is an isomorphism of\linebreak
$\textrm{Lip}_H^{\alpha}(G, \,d)$ onto the closed subalgebra of
$\textrm{Lip}^{\alpha}(M(A), \,\sigma)$ consisting of those
functions in $\textrm{Lip}^{\alpha}(M(A),\, \sigma)$ that are
continuous on $M(A)$ and are holomorphic on $\textrm{Int}\,M(A)$
w.r.t. the Gel'fand topology. In particular, for each
$\a\,\in\,(0,\,1]$, $\textrm{Lip}_H^\a(G,\,d)$ is a commutative
semisimple Banach algebra.
\end{thm}
\noindent {\bf Proof.} If $f\,\in\,\textrm{Lip}_H^\a(G, \,d)$ and
$\|f\|_{\textrm{Lip}^\a}\,\leq\,1$, then $\|f\|_\a\,\leq\,1$,
so\linebreak $|f(z) - f(w)|\,\leq\,d(z,\,w)^\a$ for all
$z,\,w\,\in\,G$. Thus $\sigma(z,\,w)\,\leq\,d(z,\,w)^\a$ for all
$z,\,w\,\in\,G$. Hence for any $f\,\in\,\textrm{Lip}_H^\a(G,\,d)$,
$\|f\|_\a\,\leq\,|\|\hat{f}\||_\a$. Since each
$\hat{f},\;f\,\in\,\textrm{Lip}_H^\a(G,\,d)$, is continuous on
$M(A)$ w.r.t. the Gel'fand topology and since $G$ is dense in
$M(A)$ w.r.t. the Gel'fand topology by Lem. 3.2, we have
$\|f\|_\infty\,=\,|\|\hat{f}\||_\infty$. Thus for all
$f\,\in\,\textrm{Lip}_H^\a(G,\,d)$,
$$\|f\|_{\textrm{Lip}^\a}\,=\,\|f\|_\infty\,+\,\|f\|_\a\,\leq\,|\|\hat{f}\||_\infty\,+\,|\|\hat{f}\||_\a\,=\,|\|\hat{f}\||_{\textrm{Lip}^\a}.$$
This together with the inequality from Prop. 2.1 of [Sh1] gives
$$\|f\|_{\textrm{Lip}^\a}\,\leq\,|\|\hat{f}\||_{\textrm{Lip}^\a}\,\leq\,2\|f\|_{\textrm{Lip}^\a}.$$ Hence the mapping
$f\,\rightarrow\,\hat{f}$ is a bicontinuous isomorphism, and the
image of $\textrm{Lip}_H^\a(G,\,d)$ is therefore a closed
subalgebra of $\textrm{Lip}^\a(M(A),\,\sigma)$.

Let $g\,\in\,\textrm{Lip}^\a(M(A),\,\sigma)$ be continuous on
$M(A)$ and holomorphic on $\textrm{Int}\,M(A)$ w.r.t. the Gel'fand
topology and let $f\, =\, g|G$ denote the restriction of $g$ to
$G$. Then $f\,\in\,\textrm{Lip}_H^\a(G,\, d)$ as $\sigma(z,\,
w)\,\leq\,d(z,\, w)^\a$ for all $z,\, w\,\in\,G$; and $\hat{f}\,
=\, g$ since both are continuous on $M(A)$ and holomorphic on
$\textrm{Int}\,M(A)$ w.r.t. the Gel'fand topology and agree on the
dense subset $G$ by Lem. 3.2. Thus those
$g\,\in\,\textrm{Lip}^\a(M(A),\, \sigma)$ which are continuous on
$M(A)$ and holomorphic on $\textrm{Int}\,M(A)$ w.r.t. the Gel'fand
topology lie in the range of the mapping $f\,\rightarrow\,\hat{f}$
from $\textrm{Lip}_H^\a(G, \,d)$. Since every $\hat{f}$ is
continuous on $M(A)$ and holomorphic on $\textrm{Int}\,M(A)$
w.r.t. the Gel'fand topology, the image of $\textrm{Lip}_H^\a(G,\,
d)$ under the mapping $f\,\rightarrow\,\hat{f}$ is exactly the set
of functions in $\textrm{Lip}^\a(M(A),\, \sigma)$ which are
continuous on $M(A)$ and holomorphic on $\textrm{Int}\,M(A)$
w.r.t. the Gel'fand topology. $\hfill \Box$

The root of automatic continuity theory is the
`uniqueness-of-norm' problem, which asks which Banach algebras
have a unique complete norm. We remark that Rickart in 1950 raised
a historically important question whether or not each semisimple
Banach algebra has a unique complete norm [Ri]. This was resolved
by Johnson in 1967 [J1]; he also established the uniqueness of the
complete norm of Banach algebras of power series in the same year
[J2]. In the more general case of Fr\'{e}chet algebras, the
uniqueness of the Fr\'{e}chet algebra topology of Fr\'{e}chet
algebras of power series was established by the author in [P1]
(see also [P3] for this result for Fr\'{e}chet algebras of power
series in several variables). In [DPR], we obtain other
interesting automatic continuity results for Fr\'{e}chet algebras
of power series, including an affirmative answer to the
Dales-McClure problem (1977) in a stronger sense, as well as a
progressive result on the most prestigious Michael problem (1952)
that the {\it test algebra} for this problem is, indeed, a
Fr\'{e}chet algebra of power series. As an easy consequence of the
above theorem, we have the following
\begin{cor} \label{uniqueness_norm} For each $\a\in(0, 1]$,
$\textrm{Lip}_H^\a(G, d)$ has a unique complete norm.$\hfill \Box$
\end{cor}

The Banach algebra $A$ is called {\it regular} if for each proper
closed subset $K$ of $M(A)$ and each point $\phi \in M(A) - K$
there exists an $f \in A$ such that $\hat{f}(\phi) = 1$ and
$\hat{f}(K) = 0$. For each $\a \in (0, 1]$, the Banach algebra
$\textrm{Lip}^\a(X, d^{'})$ is regular follows from the following
proposition, due to J. Lindberg (unpublished); for proof, see
[Sh1]. An algebra $A$ of functions on a space $X$ is {\it closed
under truncation} if the function $\min(f,\, 1)\,\in\, A$ for all
real-valued $f\,\in\,A$.
\begin{prop} \label{Lindberg_regular} Let $X$ be a compact Hausdorff space and $A$ be
a self-adjoint subalgebra of $C(X)$ which separates the points of
$X$ and contains the constant functions. If $A$ is closed under
truncation, then $A$ is regular.
\end{prop}
As a consequence, we have the following
\begin{cor} \label{Lipschitz_regular} For each $\a\,\in\,(0,\,1]$, the
Banach algebra $A\,=\,\textrm{Lip}^\a(X,\,d^{'})$ is regular.
\end{cor}
\noindent {\bf Proof.} For each $\a \in (0, 1]$, let $f \in
\textrm{Lip}^\a(X, d^{'})$ be a real-valued, distance function
defined on $X$ for a fixed point $t \in X$ and set $Tf\,=\,\min(f,
\,1)$. Then $|(Tf)(x) - (TF)(y)|\,\leq\,|f(x) - f(y)|$ for all $x,
\,y\,\in\,X$ which may be seen by comparing the graphs of $f$ and
$Tf$, or by checking each of the possible cases for a given $x$
and $y$. It is immediate from this that $\textrm{Lip}^\a(X,\,
d^{'})$ is closed under truncation. So $\textrm{Lip}^\a(X,\,
d^{'})$ is a point-separating algebra of functions on $X$. Hence
each $x\,\in\,X$ can be identified with the point evaluation
functional $\phi_x$ in the maximal ideal space $M(A)$ where
$\phi_x(f)\,=\,f(x)$, and so, the injection mapping
$x\,\ra\,\phi_x$ is one-to-one from $X$ to $M(A)$ and we may
regard $X$ as a subset of $M(A)$. Clearly, it contains the
constants, and is self-adjoint. Moreover, it is easy to show that
it is an inverse-closed algebra. By the general theory of function
algebras, $X$ is dense in $M(A)$ w.r.t. the Gel'fand topology. It
follows that $\textrm{Lip}^\a(X,\, d^{'})$ is closed under
truncation if and only if
$\{\hat{f}\,:\,f\,\in\,\textrm{Lip}^\a(X,\,d^{'})\}$ is closed
under truncation. Hence, for $\a\,\in\,(0,\,1]$, $\textrm{Lip}^\a
(X,\, d^{'})$ is regular. A little reflection of all this facts
also reveals that the relative Gel'fand topology of $X$ and the
original $d-$topology of $X$ are equivalent as well. $\hfill \Box$

We remark that the Banach algebras $\textrm{Lip}_H^{\alpha}(G,
\,d)$ are certainly not regular since they are not self-adjoint.
In fact, this follows from the following converse of Prop. 3.6 and
the method of proof of the well-known real and complex
Stone-Weierstrass Theorems.
\begin{prop} \label{converse_2.4} Let $X$ be a compact Hausdorff
space and $A$ be a closed subalgebra of $C(X)$ which is closed
under truncation and contains the constant functions. If $A$ is
regular, then $A$ separates the points of $X$ and is self-adjoint.
$\hfill \Box$
\end{prop}

As a subset of $M(A)$, $G$ inherits the Gel'fand and the metric
topologies of $M(A)$. We now compare these inherited topologies of
$G$ with its original $d$-topology. We first see that the relative
Gel'fand topology of $G$ and the $d$-topology of $G$ are
equivalent by Lem. 3.2 and Prop. 3.3 above.

We now compare the two metric topologies of $G$. The next few
propositions are concerned with the relation between $d$ and
$\sigma$ on $G$; we omit the proofs in our case as the arguments
are almost same (see [Sh1, \S 3] for notions, definitions, their
interrelations and equivalence properties; also, the results on
the homomorphisms and automorphisms, obtained in \S 5 (Thm. 5.1
and Cor. 5.2), hold true in our cases: (1) for $\a\,\in\,(0,\,1]$,
$\textrm{Lip}^\a (X,\, d^{'})$, where $X$ is compact Hausdorff
space; and (2) for $\a\,\in\,(0,\,1]$, $\textrm{Lip}_H^\a (G,\,
d)$, where $G$ is relatively compact in $\cc^n$).

Since $M(A)$ lies on the unit sphere of the dual space of
$A\,=\,\textrm{Lip}_H^{\alpha}(G,\, d)$ (which is, indeed, the
space $\Lambda^*(\alpha, \infty, \infty)$ by [Ch] when $G$ is a
bounded symmetric domain in $\cc^n$), the diameter of ($M(A),\,
\sigma)$ is at most two. Thus $\sigma$ is always a bounded metric.
Moreover, since $G$ is a bounded domain in $\cc^n$, the diameter
of $(G,\, d)$ is finite and we have the following
\begin{prop} \label{metrics_bddlyequivalent} The metric $\sigma$ on $G$ is boundedly equivalent to $d$. $\hfill \Box$
\end{prop}
As a consequence, we have the following
\begin{cor} \label{metrics_bddlyequivalent1} Let $d_1$ and $d_2$
be bounded metrics on $G$. Then $A_1$ \linebreak $=\,
\textrm{Lip}_H^{\alpha}(G,\, d_1)$ and
$A_2\,=\,\textrm{Lip}_H^{\alpha}(G,\, d_2)$ have the same elements
if and only if $d_1$ and $d_2$ are boundedly equivalent. $\hfill
\Box$
\end{cor}
Moreover, we have the following
\begin{prop} \label{metrics_bddlyequivalentextn} Given the metric
space $(G, \,d)$, the Banach algebras\linebreak
$\textrm{Lip}_H^{\alpha}(G,\, d)$ and $\textrm{Lip}_H^{\alpha}(G,
\,d^{'})$, where $d^{'}\, =\,\frac{d}{1+d}$, have the same
elements and their norms are equivalent. $\hfill \Box$
\end{prop}
As a consequence, we have the following
\begin{cor} \label{metrics_equivalent} The metrics $d$ and
$\sigma$ on $G$ are always uniformly equivalent.$\hfill \Box$
\end{cor}
As an application of the above results, we have the following
\begin{thm} \label{Lip_Gleason} Let $G$ be either a smooth (strictly)
pseudoconvex domain in $\cc^n$, or a pseudoconvex domain with a
Stein neighbourhood basis, or a polydomain in $\cc^n$. Then $G$ is
Gleason solvable for Banach algebras $\textrm{Lip}_H^\alpha(G,\,
d)$, where $\alpha\,\in\,(0,\, 1]$.
\end{thm}
\noindent {\bf Proof.}  Clearly, for each $\alpha\,\in\,(0, 1]$,
$\textrm{Lip}_H^\alpha(G, d)$ is a Banach subspace of $A^2(G)$
(see Lem. 2.1 above). Then, by Prop. 2.2, there exist
$f_1,\,\dots,\,f_n$ in $A(G)$ so that $f\,=\,\sum_{j=1}^{n}(z_j -
\alpha_j)f_j$. In fact, $f_1,\,\dots,\,f_n$ are in
$\textrm{Lip}_H^\alpha(G,\, d)$, because they satisfy the
Lipschitz condition and are uniformly continuous on $(G,\, d)$;
also, $\textrm{Lip}_H^\alpha(G, \,d)$, where $\alpha\,\in\,(0,\,
1]$, is a proper subalgebra of $A(G)$ containing the identity (so,
it cannot be an ideal of $A(G)$).

Alternatively, we see that for each $\a\,\in\,(0,\,1]$,
$\textrm{Lip}_H^\a(G,\,d)$ is a semisimple, locally Stein-Banach
algebra under the usual $\|\cdot\|_{\textrm{Lip}^\a}$ norm, as
discussed in the paragraph preceding to Lem. 3.2. Hence, by Thm.
4.3 of [P4] (when applied to $\textrm{Lip}_H^\a(G,\,d)$), we see
that the Gleason problem is solvable, because its maximal ideal
space is $\overline{G}$ by Lem. 3.2. $\hfill \Box$

\noindent {\bf Remarks D.} 1. As an application of the solvability
of the corona problem (see Lem. 3.2), we have here established the
Gleason solvability of $G$ for the holomorphic H\"{o}lder and
Lipschitz spaces, answering affirmatively the question left open
by Forn\ae ss and \O vrelid in [FO] for the pseudoconvex domains
of finite type in $\cc^2$ (and in particular, such domains with a
real analytic boundary).

\noindent 2. We studied the existence of local analytic structure
in the spectrum of a Fr\'{e}chet algebra, and as a welcome bonus,
we characterized locally Riemann algebras in [P2], giving a road
map to a characterization of locally Stein algebras in order to
obtain the Gleason solvability for these algebras in [P4].

\noindent {\bf 3.2} Having solved the Gleason's problem for
holomorphic H\"{o}der and Lipschitz spaces, we next use this
result to further solve this problem for other Banach spaces of
holomorphic functions such as Lipschitz spaces of holomorphic
functions over bounded symmetric domains in $\cc^n$ [Ch],
holomorphic mean Lipschitz spaces and Hardy-Sobolev spaces on the
unit ball in $\cc^n$ [CZ], and some special classes of mixed-norm
(Bergman) spaces (including weighted Bergman spaces) of
holomorphic functions on the unit ball in $\cc^n$ [JP]. We remark
that in the net of several Banach spaces of holomorphic functions,
analysts gave more than one name to one such Banach space of
holomorphic functions; for example, what Cho and Zhu called
special classes of mixed-norm spaces (e.g., holomorphic (mean)
Lipschitz spaces) are exactly same as special classes of the
mixed-norm Bergman spaces (e.g., holomorphic (mean)
Besov-Lipschitz spaces, named by Jevti\'{c} and Pavlovi\'{c}
[JP]), which is clear from various characterizations obtained for
these spaces in [CZ, CKK, JP].

Recall that Carleson in [C] used the interpolation (by $H^\infty-$
functions) method to solve the corona problem for the open unit
disc in $\cc$. Dronov and Kaplitskii in [DK] used a combination of
the classical dead-end space method and the interpolation (by
bounded operators on the cones in certain Banach spaces) method to
solve the Mityagin's basis problem from the theory of nuclear
Fr\'{e}chet spaces. We will use this combination, but for the
interpolation (by certain Banach spaces of bounded holomorphic
functions) method to establish the Gleason solvability and the
corona theorem below. Indeed, we use the classical dead-end space
method to show the Gleason solvability for these spaces of various
(poly)domains in $\cc^n$. The central idea of this method is as
follows: we construct the triple $(F_0, \,F_\infty,\, F)$ of
Banach spaces of holomorphic functions of a certain (poly)domain
$G$ in $\cc^n$ so that $G$ is Gleason solvable for $F_0$ and
$F_\infty$ such that $F_\infty\,\subset\,F\,\subset\,F_0$
continuously and $F_0\,\cong\,F_\infty$. Furthermore, taking
appropriate Gleason solvable Banach spaces of holomorphic
functions of $G$, we can always achieve $F_0\,\subset\,E_0$
continuously, $F_\infty\,\subset\,E_\infty$ continuously,
$E_\infty\,\subset\,E_r\,\subset\,E_0$ continuously and for each
$r$ the triple $(E_0,\, E_\infty,\, E_r)$ is an interpolation
triple of Gleason solvable Banach spaces of holomorphic functions
of $G$. The system of Banach spaces thus constructed, is always a
``complete minimal" system, squeezing $F$. Thus, $F$ is Gleason
solvable, provided that the corresponding $F_0$ and $F_\infty$ are
Gleason solvable. We remark that the classical dead-end space
method has been used to affirmatively solve certain well-known,
long-standing analytical problems from (functional) analysis such
as the Mityagin's basis problem (Pelczy\'{n}ski's basis problem is
still open, and the Grothendieck's basis problem is the most
prestigious, general one among all basis problems in the theory of
nuclear Fr\'{e}chet spaces) [DK].

We first remark that for $\alpha\,\in\,(0, 1]$ and $G$ being a
bounded symmetric domain in $\cc^n$, the Banach algebras
$\textrm{Lip}_H^\alpha(G,\,d)$, discussed above, are the classical
holomorphic Lipschitz spaces
$\Lambda_\alpha\,=\,H\Lambda(\alpha/n,\,\infty,\,\infty)$,
discussed by Chen [Ch, Thm. 10] that are actually same as the
spaces $\Lambda_\alpha^{\infty,\,\infty}$, discussed by Cho and
Zhu [CZ] for the unit ball (note that we here include the case of
the Lipschitz algebra $\textrm{Lip}_H(G,\,d)$ of holomorphic
functions, which is the space
$\Lambda_1\,=\,\Lambda_1^{\infty,\,\infty}$, discussed by Cho and
Zhu [CZ], that is actually same as the space
$H\Lambda(1/n,\,\infty,\,\infty)$, discussed by Chen [Ch, Thms. 4
and 10]). Krantz named $\Lambda_\alpha$ as holomorphic Zygmund
spaces in [K1-2, KL1-2]. Moreover, by Lem. 2 of [Ch], we see that
a function $f\,\in\,H\Lambda(\alpha/n,\,\infty,\,\infty)$ is
continuous on $G\,\bigcup\,b$, where $b$ is a Bergman-\u{S}ilov
boundary of $G$ (i.e., $f$ is uniformly continuous on $G$ by [A]),
and by Thm. 6 of [Ch], for $1\,\leq\,p\,<\,\infty$ and
$q\,=\,\infty$,
$\Lambda_\alpha\,\subset\,H\Lambda(\alpha,\,p,\,\infty)\,\subset\,\Lambda_{\alpha
- 1/p}$ (and these inclusions are continuous). Since
$\Lambda_\alpha\,\cong\,\Lambda_{\alpha - 1/p}$ by Thm. 4 of [Ch],
we have that $G$ is Gleason solvable for
$H\Lambda(\alpha,\,p,\,\infty)$ by the classical dead-end space
method, discussed above (note that one suitably chooses $E_0,
\,E_\infty$ and $E_r$ to be the classical holomorphic Lipschitz
spaces (which are Gleason solvable), because for
$0\,<\,\alpha\,<\,\beta\,\leq\,1$, we have
$\Lambda_\beta\,\subset\,\Lambda_\alpha$ continuously). Along the
same lines, $G$ is Gleason solvable for
$H\Lambda(\alpha,\,\infty,\,q)$ for $1\,\leq\,q\,<\,\infty$,
since, by Thm. 5 of [Ch], for $\alpha\,<\,1$, $p\,=\,\infty$ and
$1\,\leq\,q\,<\,s\,<\,\infty$,
$\Lambda_1\,\subset\,H\Lambda(\alpha,\,\infty,\,q)\,\subset\,H\Lambda(\alpha,\,\infty,\,s)\,\subset\,\Lambda_\alpha$
(and these inclusions are continuous).

Concerning the Gleason solvability for holomorphic mean Lipschitz
spaces, discussed by Cho and Zhu in [CZ], we first remark that for
$0\,<\,\alpha_1\,<\,\alpha_2\,<\,\infty$,
$1\,\leq\,p,\,q_1\,<\,\infty$ and $q_2\,=\,\infty$, we have
\begin{equation} \label{Prop 3.6} \Lambda_{\alpha_2}^{p,\,q_1}\,\subset\,\Lambda_{\alpha_2}^p\,\subset\,\Lambda_{\alpha_1}^{p,\,q_1}\,\subset\,\Lambda_{\alpha_1}^p,
\end{equation}
by (i) and (ii) of Prop. 3.6 of [CZ] (and all this inclusions are
continuous; also, recall the roles of the indices $\alpha,\, p$
and $q$ in the holomorphic mean Lipschitz space
$\Lambda_\alpha^{p,\,q}$ as follows: The index $p$ indicates the
basic $H^p$ norm that is used, the index $\a$ gives the order of
smoothness involved, and the index $q$ represents a rather subtle
correction to the order of smoothness; Prop. 3.6 of [CZ] makes
these remarks somewhat more precise). In fact, by Thm. 3.12 of
[CZ], all this spaces are isomorphic. Further, it is easy to see
that Thm. 3.12 of [CZ] also holds true for the case
$\alpha\,>\,0$, $s\,=\,q\,=\,\infty$ and $\beta\,=\,1/p$ with
$1\,\leq\,p\,<\,\infty$; i.e.,
$\Lambda_\alpha\,\cong\,\Lambda_{\alpha - 1/p}$. In Thm. 6 of
[Ch], Chen discussed certain inclusion relations between the
classical holomorphic Lipschitz spaces and holomorphic, integrated
Lipschitz spaces (i.e., $p\,<\,\infty$). The situation for the
classical holomorphic mean (integrated) Lipschitz spaces will be
similar, and thus, we have
\begin{equation} \label{Thm 6} \Lambda_\alpha\,\subset\,\Lambda_\alpha^p\,\subset\,\Lambda_{\alpha - 1/p}
\end{equation} (and these inclusions are continuous). So, the space $\Lambda_\alpha^p$ is squeezed
between two classical holomorphic mean Lipschitz spaces. Summing
up all above situations, we see that the unit ball is Gleason
solvable for all holomorphic mean (integrated) Lipschitz spaces
(this includes the special case of the diagonal holomorphic Besov
spaces $B_p\,=\,\Lambda_{\frac{n}{p}}^{p, p}$ with $p\,<\,\infty$
and $\frac{n}{p}\,<\,1$).

Concerning the Gleason solvability for Hardy-Sobolev spaces,
discussed by Cho and Zhu in [CZ], we have two cases: (1)
$p\,\in\,(1, 2]$, and (2) $p\,\in\,[2, \infty)$. For the case (1),
by (i) of Prop. 3.6 and (c) and (b) of Thm. 5.2 of [CZ], we have
$\Lambda_\alpha^{p, p}\,\subset\,
H_\alpha^p\,\subset\,\Lambda_\alpha^{p,
2}\,\subset\,\Lambda_\alpha^p$. For the case (2), by (i) of Prop.
3.6 and (a) and (d) of Thm. 5.2 of [CZ], we have
$\Lambda_\alpha^{p, 2}\,\subset\,
H_\alpha^p\,\subset\,\Lambda_\alpha^{p,
p}\,\subset\,\Lambda_\alpha^p$. So, the unit ball is also Gleason
solvable for all Hardy-Sobolev spaces.

Next, we obtain the Gleason solvability for some special classes
of holomorphic (mean) Besov-Lipschitz spaces, discussed by
Jevti\'{c} and Pavlovi\'{c} in [JP]. For this, we note that for
$0\,<\,\alpha\,<\,n$ and $\phi$ an increasing, positive,
continuous function on $(0, \,1]$, the space $\Lambda_{n, \phi}^p$
contains the Hardy-Sobolev space $H_n^p$, and the space
$\Lambda_{n, \phi}^{p, q}$ contains the Besov-Hardy-Sobolev space
$H_n^{p, q}$. Further, if the function $x\,\mapsto
\phi(x)/x^\alpha$ increases on $(0, \,1]$, then $\Lambda_{n,
\phi}^p\,\subset\,\Lambda_\alpha^p$; in fact, we have $\Lambda_{n,
\phi}^p\,\cong\,\textrm{Lip}_{n, \phi}^p$, by Cor. 1.9, (1.12) and
(1.16) of [JP]. The same situations hold true for the integrated
spaces $\Lambda_{n, \phi}^{p, q}$, $\Lambda_\alpha^{p, q}$ and
$\textrm{Lip}_{n, \phi}^{p, q}$, by Cor. 1.10, (1.11), (1.18) and
(1.19) of [JP]. Since the spaces $\Lambda_{n, \phi}^p$ and
$\Lambda_{n, \phi}^{p, q}$ are squeezed between
(Besove)-Hardy-Sobolev spaces and holomorphic mean
(Besov-)Lipschitz spaces (which are Gleason solvable), we see that
the unit ball is Gleason solvable for these special classes of
holomorphic (mean) Besov-Lipschitz spaces.

Finally, it is easy to obtain analogues of (4.1), (4.2) and Thm. 4
of [Ch] for Hardy-Sobolev spaces and Besov-Hardy-Sobolev spaces to
show that the unit ball is Gleason solvable for
Besov-Hardy-Sobolev spaces (including weighted Bergman spaces)
since so are Hardy-Sobolev spaces.

Thus, summarizing the above discussion, we have the following
\begin{thm} \label{Gleason_Solvability} Let $G$ be a bounded
domain in $\cc^n$. Then $G$ is Gleason solvable for the following
Banach spaces of holomorphic functions on $G$.
\begin{enumerate}
\item[{\rm (i)}] holomorphic (Besov-)Lipschitz spaces on a bounded
symmetric domain $G$ in $\cc^n$; \item[{\rm (ii)}] holomorphic
mean (Besov-)Lipschitz spaces and Hardy-Sobolev spaces on the unit
ball in $\cc^n$; \item[{\rm (iii)}] some special classes of
holomorphic mixed-norm (Bergman) spaces (including weighted
Bergman spaces) on the unit ball in $\cc^n$; and \item[{\rm (iv)}]
some special classes of holomorphic (mean) Besov-Lipschitz spaces
such as the spaces $\Lambda_{n, \phi}^p$ and the integrated spaces
$\Lambda_{n, \phi}^{p, q}$, where a function $\phi$ is as above,
on the unit ball in $\cc^n$. $\hfill \Box$
\end{enumerate}
\end{thm}
\section{Corona problem for Banach spaces of bounded holomorphic functions.}
In this section, we solve the corona problem for various Banach
spaces of bounded holomorphic functions, as we have done in the
previous section for the Gleason solvability. We give these
solutions by various (functional analytic) methods, and for a
number of suitable bounded domains in $\cc^n$, including the ball
and the polydisc. We recall that we have already discussed
solution of the original corona problem for the Banach algebra
$H^\infty(G)$, where $G$ is a certain pseudoconvex domain in
$\cc^n$ and a polydomain in $\cc^n$ (including the open unit ball
and polydisc), in Theorem 2.4 (and Corollary 2.5) above; we use
the Gleason solvability and a Banach algebraic technique to obtain
this solution.

We note that Krantz and Li restated the corona problem in the
language of a certain bounded operator which needs to be onto, if
one wants to solve the original problem [KL1], and with this
approach, they succeed to solve the corona problem for Bergman
spaces $A^p(G), \;0\,<\,p\,<\,\infty$, where $G$ is the unit ball
as well as a bounded strictly pseudoconvex domain in $\cc^n$ with
$C^3-$boundary (but for the range $1\,<\,p\,<\,\infty$), for Hardy
spaces $H^p(G), \;0\,<\,p\,<\,\infty$, where $G$ is a bounded
strictly pseudoconvex domain in $\cc^n$ with $C^3-$boundary, and
for holomorphic Zygmund spaces
$H\Lambda_\a(G),\;0\,<\,\a\,<\,\infty$, where $G$ is a bounded
strictly pseudoconvex domain in $\cc^n$ with $C^3-$boundary (see
Thms. 1.1 - 1.4 in [KL1]). We remark that the corona data
$f\,=\,(f_1,\,\dots,\,f_m)$ belongs to $H^\infty(G)$ in Thms. 1.1
- 1.3, but it is easy to see that they could have taken the corona
data in the respective Banach spaces as well since for
$1\,\leq\,p\,<\,\infty$, $H^\infty\,\hookrightarrow\,H^p$
continuously. We also remark that they could not solve the
original corona problem (see Remark 1 of [KL1]); also, they used
the $\overline{\partial}-$technique to prove these theorems, which
is not always possible for arbitrary strictly pseudoconvex domains
(see, also, [V]). Indeed, we show that the bounded operator $S_f$
is onto, solving the original corona problem in the following
\begin{thm} \label{KL_Sf} Let $G$ be a strictly pseudoconvex domain
in $\cc^n$ with \linebreak $C^3-$boundary or a polydomain in
$\cc^n$. Let $f\,=\,(f_1,\,\dots,\,f_m)$ be the corona data from
the product space $H^\infty(G)_m$. Then
$S_f\,:\,H^\infty(G)_m\,\rightarrow\,H^\infty(G)$ is onto. In
particular, the corona problem holds true in $H^\infty(G)$.
\end{thm}
\noindent {\bf Proof.} This proof is, indeed, a little reflection
of the proof of Thm. 2.3 above. First, we note that $H^\infty(G)$
is a semi-simple Banach algebra of power series in $\F_n$,
satisfying the hypotheses of Lem. 3.1 of [P4], because $G$ is
Gleason solvable by Thm. 2.3 above, and
$\bigcap_{k\,\geq\,1}\overline{M^k}\,=\,\{0\}$. Now, following
proof of (iii) of Lem. 3.1 of [P4], we see that the product space
$H^\infty(G)_m$ is, indeed, a Banach algebra with an identity; the
algebraic operations are coordinatewise. Further, with a minor
modification, we also see that for each
$g\,=\,(g_1,\,\dots,\,g_m)\,\in\,H^\infty(G)_m$,
$T_g\,:=\,\oplus_{i =
1}^mT_{g_i},\;\oplus_{i\,=\,1}^mh_i\,\mapsto\,\sum_{i\,=\,1}^m(f_i
- g_i)h_i,$ is a bounded linear map on $H^\infty(G)_m$, where
$T_{g_i}\,:\,h\,\mapsto \,(f_i - g_i)h$ is a bounded linear
operator on $H^\infty(G)$. Now, given the corona data
$f\,=\,(f_1,\,\dots,\,f_m)$, let $I$ be the algebraically finitely
generated ideal in $H^\infty(G)$. Then $I$ is an improper ideal if
and only if $I$ cannot be contained in any maximal ideal of
$H^\infty(G)$, and $\sum_{i\,=\,1}^mf_ih_i\,=\,1$ for some
$h\,=\,(h_1,\,\dots,\,h_m)\,\in\,H^\infty(G)_m$. That is, $S_f$ is
onto (clearly, for the corona data $f\,=\,(f_1,\,\dots,\,f_m)$,
$S_f$ is obviously a bounded linear map; in fact, $S_f \,=\,T_0$).
$\hfill \Box$

Next, we recall that we have already discussed solution of the
corona problem for semi-simple, locally Stein-Banach algebras
$\textrm{Lip}_H^\a(G,\, d)$, where $\a\,\in\,(0, 1]$ and $G$ is a
relatively compact, ``pseudoconvex" domain in $\cc^n$ as in Lem.
3.2 (including a relatively compact symmetric domain, and in
particular, the unit ball). In fact, in this case, we have better
result in the sense that the relative Gel'fand topology of $G$ is
equivalent to the $d-$ and $\sigma-$ topologies of $G$.
Alternatively, one applies the Banach-algebraic technique,
discussed in Thm. 2.4 or, chooses to show that $S_f$ is onto as
above, to establish the corona theorem for these algebras. In
particular, we have solved the question of Forn\ae ss and \O
vrelid [FO] for the Gleason solvability for pseudoconvex domains
of finite type in $\cc^2$ (and in particular, such domains with a
real-analytic boundary); see Thm. 3.13 above. Also, by Lem. 3.2,
the corona theorem holds true in $\textrm{Lip}_H^\a(G,\, d)$,
where $\a\,\in\,(0, 1]$ and $G$ is the ball or, polydisc in
$\cc^n$. This would extend the result of [KL2] for the polydisc,
giving us a non-linear solution (see their remark before
References in [KL2]).

Now, using the classical dead-end space method, discussed in \S
3.2, we establish the corona theorem for various Banach spaces of
bounded holomorphic functions on bounded symmetric domains in
$\cc^n$ (including the unit ball). We do not repeat the arguments
for these spaces, but remark that if $F$ is the middle space in
the inclusion relation for which we wish to establish the corona
theorem for a given corona data in that space, then it is clear
that we have the same corona data in $F_0$ (and hence in
$F_\infty$, too, as we have $F_0\,\cong\,F_\infty$). Now, since
the corona theorem holds true in both $F_\infty$ and $F_0$, it
also holds true in $F$ as well. Thus, we have the following
\begin{thm} \label{Corona_Solvability} Let $G$ be a bounded
domain in $\cc^n$. Then the corona theorem holds true for $G$ for
the following Banach spaces of bounded holomorphic functions on
$G$.
\begin{enumerate}
\item[{\rm (i)}] (Besov-)Lipschitz spaces of holomorphic functions
over a bounded symmetric domain $G$ in $\cc^n$; \item[{\rm (ii)}]
holomorphic mean (Besov-)Lipschitz spaces and Hardy-Sobolev spaces
on the unit ball in $\cc^n$; \item[{\rm (iii)}] some special
classes of holomorphic mixed-norm (Bergman) spaces (including
weighted Bergman spaces) on the unit ball in $\cc^n$; and
\item[{\rm (iv)}] some special classes of holomorphic (mean)
Besov-Lipschitz spaces such as the spaces $\Lambda_{n, \phi}^p$
and the integrated spaces $\Lambda_{n, \phi}^{p, q}$, where a
function $\phi$ is as above, on the unit ball in $\cc^n$. $\hfill
\Box$
\end{enumerate}
\end{thm}
\section{Concluding remarks on our approaches and open problems.}
As we explain in \S 1, given a finite corona date, if one wants to
solve the corona problem in the SCV, one has to solve the
B\'{e}zout's equation (the algebraic formulation), and analysts
mostly considered the approach of solving the
$\bar{\partial}-$equation (the classical analysis method), but
this technique was failed in the past while solving the original
corona problem for $H^\infty(G)$, even when $G$ is the open unit
ball in $\cc^n$ [V]. Then, if one wants to solve this problem by
considering its equivalent formulation (the topological
formulation), one may consider the approach of solving the {\bf
ideal problem} (which is of algebraic nature). Thus this approach
clearly indicates/hints the possible connection between the corona
problem and the Gleason's problem in the SCV case. However, until
we consider this connection for the first time in Thm. 4.3 of
[P4], as far as we know, no analysts did take/consider this
approach, and it is not a great surprise for this for obvious
reasons. Then, our extension of the Clos' result also plays a
significant role while filling the gap via the Gleason solvability
of certain domains $G$ for $H^\infty(G)$ when one wants to take
the approach of solving the ideal problem by using the
Banach-algebraic technique (which actually provides the necessary
and sufficient conditions for the ideal problem [DKSTW, \S 7]).
Thus, we consider a very strong connection between the two
significant problems from the theory of SCV by discussing an
equivalence between these two problems.

At this point, one may think why a very short and the most elegant
approach of Krantz and Li [KL1] could not find the solution of the
original corona problem while applying the classical analysis
method of solving the $\bar{\partial}-$equation. Again, an answer
lies in the fact that we fill the gap through the method of proof
of Lem. 3.1 of [P4] (where one needs the Gleason solvability as
one of the two hypotheses) and the Banach-algebraic technique.

Finally, we give some comments on our approach of the combination
of the Banach algebra method and the classical dead-end space
method to solve both the problems in the SCV for various familiar
Banach spaces of (bounded) holomorphic functions on certain
domains in $\cc^n$. First, we note that due to the algebraic
nature of both the problems, it is very much clear that one does
not require the Banach algebra structure of the underlying Banach
space of (bounded) holomorphic functions for which we want to
solve both the problems; equivalently, analysts did give solutions
of both the problems for certain Banach spaces of certain domains
by applying various classical analysis methods. However, as we
remark above, those classical methods are hard and complicated in
the SCV case. However, in \S 3, we first establish both the
problems for the Lipschitz algebras $\textrm{Lip}_H^\a(G,\,d)$ of
bounded holomorphic function on certain domains $G$ in $\cc^n$,
where $\a\,\in\,(0,\,1]$. This provides a concrete base for the
classical dead-end space method, in order to discuss solutions of
both the problems for various familiar Banach spaces of (bounded)
holomorphic functions on certain domains in $\cc^n$ (and this set
of Banach spaces are complemented in the sense that we mostly do
not repeat Banach spaces for which solutions of both the problems
are available, as mentioned above).

\noindent {\bf Open problems.} We now discuss some future
directions/open problems in view of our approaches as follows.
\begin{enumerate}
\item[{\rm (1)}] Throughout the paper, we have considered certain
domains $G$ from $\cc^n$ only. As commented in \S 1, our abstract
(functional analytic) method from [P4] has been applied to deduce
the Gleason $A$-property of a semisimple, locally Stein-Banach
algebra in Thm. 4.3 or, Cor. 4.4, where $Y$ may be a domain (i.e.,
a (reduced) Stein space) from the complex manifolds (of finite
dimension). Moreover, the Banach algebra method is also available
for a more general situation. Hence we wish to solve the corona
problem for such domains for various locally Stein-Banach algebras
of bounded holomorphic functions defined on these domains.

\item[{\rm (2)}] In 4.4 of [P4], we have raised the Gleason's
problem (as well as the Gleason solvability problem) in the
infinite dimensional case for the Banach algebras $H^\infty(B_X)$
and ${\cal A}_u(B_X)$ (a subalgebra of $H^\infty(B_X)$), where
$B_X$ is the open unit ball of a complex Banach space $X$. It is
interesting to further develop our approach (i.e., the abstract
(functional analytic) method and an equivalency of the two
problems so as to solve both the problems for the two Banach
algebras.

\item[{\rm (3)}] Is the Hartogs triangle $T$ Gleason solvable for
$H^\infty(T)$ (or, $A(T)$)? Similarly, does the corona theorem
hold true for $H^\infty(G)$, where $G$ is either the worm domain
or, the Hartogs triangle? Recall that Hartogs triangle does not
have a smooth boundary, and the worm domain does not have a Stein
neighbourhood basis [DF].
\end{enumerate}

\noindent Address: Department of Mathematics, C. U. Shah
University, Wadhwan City, Gujarat, INDIA.

\noindent E-mails: srpatel.math@gmail.com, coolpatel1@yahoo.com

\begin{thebibliography}{Dillo 83}
\bibitem[A]{A} R. M. Aron, \emph{Algebras of bounded analytic functions on finite and infinite dimensional balls}, Lecture series, workshop on Advanced Functional Analysis and its Applications, IIT Hyderabad, 2020.

\bibitem[ACG]{ACG} R. M. Aron, B. J. Cole and T. W. Gamelin, \emph{Spectra of algebras of analytic functions on a Banach space}, J. Reine Angew. Math. 415 (1991), 51-93.

\bibitem[ADLM]{ADLM} R. M. Aron, V. Dimant, S. Lassale and M. Maestre, \emph{Gleason parts for algebras of holomorphic functions in infinite dimensions}, Rev. Mat. Complut. 33 (2020), 415-436.

\bibitem[AFGM]{AFGM} R. M. Aron, J. Falc\'{o}, D. Garc\'{i}a and M. Maestre, \emph{Analytic structure in fibers}, Studia Math. 240 (2018), 101-121.

\bibitem[BF]{BF} U. Backlund and A. F\"{a}llstr\"{o}m, \emph{The Gleason property for Reinhardt domains}, Math. Ann. 308 (1997), 85-91.

\bibitem[Be]{Be} F. Beatrous Jr, \emph{H\"{o}lder estimates for the $\bar{\partial}$-equation with a support condition}, Pacific J. Math. 90 (1980), 249-257.

\bibitem[BO]{BO} A. Beddaa and M. Oudadess, \emph{On a question of A. Wilansky in normed algebras}, Studia Math. T. XCV (1989), 175-177.

\bibitem[B]{B} B. Berndtsson, \emph{A smooth pseudoconvex domain in $\cc^2$ for which $L^\infty-$estimates for $\bar{\partial}$ do not hold}, Ark. Mat. 31 (1993), 209-218.

\bibitem[BCL] {BCL} B. Berndtsson, S. Y. Chang and K. C. Lin, \emph{Interpolating sequences in the polydisc}, Trans. Amer. Math. Soc. 302 (1987), 161-169.

\bibitem[Br]{Br} A. Brudnyi, \emph{Corona problem for $H^\infty$ on Riemann surfaces}, Fields Institute Communications 72 (2014), 31-45.

\bibitem[C]{C} L. Carleson, \emph{Interpolation by bounded analytic functions and the corona problem}, Ann. Math. (2) 76 (1962), 547-559.

\bibitem[Ca]{Ca} R. L. Carpenter, \emph{Principal ideals in F-algebras}, Pacific J. Math. 35 (1970), 559-563.

\bibitem[Ch]{Ch} W-Y. Chen, \emph{Lipschitz spaces of holomorphic functions over bounded symmetric domains in $\cc^n$}, J. Math. Anal. Appl. 81 (1981), 63-87.

\bibitem[CZ]{CZ} H. R. Cho, and K. Zhu, \emph{Holomorphic mean Lipschitz spaces and Hardy Sobolev spaces on the unit ball}, Complex Variables and Elliptic Equations 57 (2012), 995-1024.

\bibitem[CKK]{CKK} H. R. Cho, H. W. Koo and E. G. Kwon, \emph{Holomorphic mean Lipschitz spaces and Besov spaces on the unit ball in $\cc^n$}, J. Korean Math. Soc.

\bibitem[Cl]{Cl} T. G. Clos, \emph{Solvability of the Gleason problem on a class of bounded pseudoconvex domains}, Complex Anal. Oper. Theory 16, no. 58 (2022), 1-9.

\bibitem[CI]{CI} T. G. Clos and A. J. Izzo, \emph{Approximation by an algebra generated by holomorphic and conjugate holomorphic functions}, Rocky Mountain J. Math. 52 (2022), 1289-1294.

\bibitem[CGJ]{CGJ} B. J. Cole, T. W. Gamelin and W. B. Johnson, \emph{Analytic disks in fibers over the unit ball of a Banach space}, Michigan Math. J. 39 (1992), 551-569.

\bibitem[DPR]{DPR} H. G. Dales, S. R. Patel and C. J. Read, \emph{Fr\'{e}chet algebras of power series}, In: Loy, R.J., Runde, V., Soltysiak, A. (eds.) Banach Algebras 2009, pp. 123-158, Banach Center Publi. 91, Polish Acad. Sci., Warsaw, 2010.

\bibitem[DF]{DF} K. Diederich and J. E. Fornæss, \emph{Pseudoconvex domains: An exmaple with non-trivial nebenh¨ulle}, Math. Ann. 225 (1977), 275-292.

\bibitem[DKSTW]{DKSTW} R. G. Douglas, S. G. Krantz, E. T. Sawyer, S. Treil and B. D. Wick (Eds.), \emph{The Corona Problem, Connections Between Operator Theory, Function Theory, and Geometry Series}, Fields Institute Communications, 72 (2014), 1-29.

\bibitem[DK]{DK} A. K. Dronov and V. M. Kaplitskii, \emph{On the existence of a basis in a complemented subspace of a nuclear K\"{o}the space in the class $(d_1)$}, (Russian) Mat. Sb. 209 (2018), 50-70; translation in Sb. Math. 209 (2018), 1463-1481.

\bibitem[F]{F} T. M. Flett, \emph{Lipschitz spaces of funtions on the circle and the disc}, J. Math. Anal. Appl. 39 (1972), 125-158.

\bibitem[FO]{FO} J. E. Forn\ae ss and N. \O vrelid, \emph{Finitely generated ideals in $A(\Omega)$}, Ann. Inst. Fourier (Grenoble) 33 (1983), 77-85.

\bibitem[FS]{FS} J. E. Forn\ae ss and N. Sibony, \emph{Smooth pseudoconvex domains in $\cc^2$ for which the corona theorem and $L^p$ estimates for $\bar{\partial}$ fail}, Complex Analysis and Geometry, 209-222, Univ. Ser. Math., Plenum, New York, 1993.

\bibitem[G]{G} J. B. Garnett, \emph{Bounded Analytic Functions}, Academic Press, New York, 1981.

\bibitem[Gl]{Gl} A. M. Gleason, \emph{Finitely generated ideals in Banach algebras} J. Math. Mech. 13 (1964), 125-132.

\bibitem[Go]{Go} H. Goldmann, \emph{Uniform Fr\'{e}chet algebras}, North Holland Publ. Co., Amsterdam, 1990.

\bibitem[HS]{HS} M. Hakim and N. Sibony, \emph{Spectre de $A(\bar{\Omega})$ pour les domains borne\'{e}s faiblement pseudoconvexes r\'{e}guliers}, J. Funct. Anal. 37 (1980), 127-135.

\bibitem[H1]{H1} G. M. Henkin, \emph{Integral representations of functions holomorphic in strictly pseudo-convex domains and some applications}, Math. USSR, Sb. 7 (1970), 597-616 (English).

\bibitem[H2]{H2} G. M. Henkin, \emph{Approximation of functions in strictly pseudoconvex domains and a theorem of Z.L. Leibenzon}, Bull. Acad. Polon. Sci S\'{e}r. Math. Astronom. Phys. 19 (1971), 37-42 (in Russian).

\bibitem[Ho]{Ho} K. Hoffman, \emph{Banach spaces of analytic functions}, Prentice-Hall, Englewood Cliffs, New Jersey, 1962.

\bibitem[Hu]{Hu} Z. Hu, \emph{Gleason's problem for harmonic mixed norm and Bloch spaces in convex domains}, Math. Nachr. 279 (2006), 164-178.

\bibitem[J]{J} P. Jak\'{o}bczak, \emph{Extension and decomposition operators in products of strictly pseudoconvex sets}, Annales Polon. Math. 44 (1984), 219-237.

\bibitem[JP]{JP} M. Jevti\'{c} and M. Pavlovi\'{c}, \emph{Besov-Lipschitz and mean Besov-Lipschitz spaces of holomorphic functions on the unit ball}, Potential Anal. 38 (2013), 1187-1206.

\bibitem[J1]{J1} B. E. Johnson, \emph{The uniqueness of the (complete) norm topology}, Bull. Amer. Math. Soc. 73 (1967), 537-539.

\bibitem[J2]{J2} B. E. Johnson, \emph{Continuity of linear operators commuting with continuous linear operators}, Trans. Amer. Math. Soc. 128 (1967), 88-102.

\bibitem[KN]{KN} N. Kerzman and A. Nagel, \emph{Finitely generated ideals in certain function algebras}, J. Funct. Anal. 7 (1971), 212-215.

\bibitem[K1]{K1} S. G. Krantz, \emph{Lipschitz spaces, smoothness of functions, and approximation theory}, Expositiones Math. 3 (1983), 193-260.

\bibitem[K2]{K2} S. G. Krantz, \emph{On a theorem of Stein}, Trans. Amer. Math. Soc. 320 (1990), 625-642.

\bibitem[K3]{K3} S. G. Krantz, \emph{Function theory of several complex variables}, AMS Chelsea Publishing, Providence, RI, 2001, Reprint of the 1992 edition.

\bibitem[K4]{K4} S. G. Krantz, \emph{Cornerstones of Geometric Function Theory: Explorations in Complex Analysis}, Birkh\"{a}user Publishing, Boston, 2006.

\bibitem[K5]{K5} S. G. Krantz, \emph{The corona problem with two pieces of data}, Proc. Amer. Math. Soc. 138 (2010), 3651-3655.

\bibitem[K6]{K6} S. G. Krantz, \emph{The corona problem in several complex variables}, Fields Institute Communications 72 (2014), 107-126.

\bibitem[KL1]{KL1} S. G. Krantz and S. Y. Li, \emph{Some remarks on the corona problem on strongly pseudoconvex domains in $\cc^n$}, Illinois J. Math. 39 (1995), 323-349.

\bibitem[KL2]{KL2} S. G. Krantz and S. Y. Li, \emph{Explicit solutions for the corona problem with Lipschitz data in the polydisc}, Pacific J. Math. 174 (1996), 443-458.

\bibitem[KL3]{KL3} S. G. Krantz and S. Y. Li, \emph{Factorization of functions in subspaces of $L^1$ and applications to the corona problem}, Indiana Univ. Math. J. 45 (1996), 83-102.

\bibitem[Li]{Li} I. Lieb, \emph{Die Cauchy-Riemannschen Differentialgleichung auf streng pseudokonveksen Gebieten: Stetige Randwerte}, Math. Ann. 199 (1972), 241-256.

\bibitem[L]{L} Y. Liu, \emph{Boundedness of the Bergman type operators on mixed norm spaces}, Proc. Amer. Math. Soc. 130 (2002), 2363-2367.

\bibitem[M]{M} H. Mirkil, \emph{Continuous translation of H¨older and Lipschitz functions}, Canad. J. Math. 12 (1960), 674-685.

\bibitem[N]{N} A. Noell, \emph{The Gleason problem for domains of finite type}, Complex Variables Theory Appl. 4 (1985), 233-241.

\bibitem[O]{O} J. M. Ortega, \emph{The Gleason problem in Bergman-Sobolev spaces}, Complex Variables Theory Appl. 20 (1992), 157-170.

\bibitem[P1]{P1} S. R. Patel, \emph{Fr\'{e}chet algebras, formal power series, and automatic continuity theory}, Studia Math. 187 (2008), 125-136.

\bibitem[P2]{P2} S. R. Patel, \emph{Fr\'{e}chet algebras, formal power series, and analytic structure}, J. Math. Anal. Appl. 394 (2012), 468-474.

\bibitem[P3]{P3} S. R. Patel, \emph{Uniqueness of the Fr\'{e}chet algebra topology on certain Fr\'{e}chet algebras}, Studia Math. 234 (2016), 31-47.

\bibitem[P4]{P4} S. R. Patel, \emph{On functional analytic approach for Gleason's problem in the theory of SCV}, Banach J. Math. Anal. 16, article no. 47 (2022), 1-17.

\bibitem[Pr]{Pr} O. Preda, \emph{The corona problem with restrictions on the relative position of sublevel sets}, Arch. Math. (Basel) 105 (2015), 563-569.

\bibitem[R]{R} E. Ramirez de Arellano, \emph{Ein Divisionsproblem und Randintegraldarstellungen in der komplexen Analysis}, Math. Ann. 184 (1970), 172-187 (German).

\bibitem[Re]{Re} T. T. Read, \emph{The powers of maximal ideals in a Banach algebra and analytic structure}, Trans. Amer. Math. Soc. 161(1971), 235-248.

\bibitem[RS1]{RS1} G. Ren and J. Shi, \emph{Bergman type operator on mixed norm spaces with applications}, Chinese Ann. Math. Ser. B 18 (1997), no. 3, 265-276, Chinese summary, Chinese Ann. Math. Ser. A 18 (1997), no. 4, 527.

\bibitem[RS2]{RS2} G. Ren and J. Shi, \emph{Gleason's problem in weighted Bergman space on egg domains}, Sci. China Ser. A 41 (1998), no. 3, 225-231.

\bibitem[Ri]{Ri} C. E. Rickart, \emph{The uniqueness of norm problem in Banach algebras}, Ann. Math. 51 (1950), 615-628.

\bibitem[Ro]{Ro} H. Rossi, \emph{Holomorphically convex sets in several complex variables}, Ann. Math 74 (1961), 470-493.

\bibitem[Ru]{Ru} W. Rudin, \emph{Function theory in the unit ball of $\cc^n$}, Classics in Mathematics, Springer-Verlag, Berlin, 2008, Reprint of the 1980 edition.

\bibitem[Sa]{Sa} S. \c{S}ahut\u{o}glu, \emph{Strong Stein neighbourhood bases}, Complex Variables and Elliptic Equations, 57 (2012), 1073-1085.

\bibitem[S1]{S1} N. Sibony, \emph{Probl\'{e}me de la couronne pour des domaines pseudoconvexes \'{a} bord lisse}, Ann. Math. (2) 126( 1987), 675-682.

\bibitem[S2]{S2} N. Sibony, \emph{Some aspects of weakly pseudoconvex domains}, Proc. Symp. Pure Math. 52 (1991), 199-231.

\bibitem[Sh1]{Sh1} D. R. Sherbert, \emph{Banach algebras of Lipschitz functions}, Pacific J. Math. 13 (1963), 1387-1399.

\bibitem[Sh2]{Sh2} D. R. Sherbert, \emph{The structure of ideals and point derivations in Banach algebras of Lipschitz functions}, Trans. Amer. Math. Soc. 111 (1964), 240-272.

\bibitem[T]{T} A. Tikaradze, \emph{On the corona problem for strongly pseudoconvex domains}, Analysis Math. 48 (2022), 1209-1212.

\bibitem[V]{V} N. Th. Varopoulos, \emph{BMO functions and the $\bar{\partial}-$equation}, Pacific J. Math. 71 (1977), 221-273.

\bibitem[X]{X} J. Xiao, \emph{The $\bar{\partial}-$problem for multipliers of the Sobolev space}, Manuscripta Math. 97 (1998), 217-232.

\bibitem[Z]{Z} K. H. Zhu, \emph{The Bergman spaces, the Bloch space, and Gleason's problem}, Trans. Amer. Math. Soc. 309 (1988), 253-268.
\end{thebibliography}
\end{document}